\documentclass[12pt]{article}
\usepackage{amsmath,enumerate,amsfonts,color,amssymb, amsthm}

\usepackage{tikz}

\def\Rot{{\mathcal R}}

\def\NN{{\mathbb N}}

\def\RR{{\mathbb R}}

\def\mcF{{\mycal F}}

\def\mcQ{{\mycal Q}}

\def\mcV{{\mycal V}}

\def\u{{\upsilon}}
\def\hot{\textrm{h.o.t.}}

\def\IdMa{{\rm Id}}

\def\ringw{{\mathring{w}}}

\def\hh{{\rm \overline H}}

\newtheorem{theorem} {\sc  Theorem\rm} [section]

\newtheorem{corollary} [theorem] {\sc  Corollary\rm}
\newtheorem{lemma} [theorem] {\sc  Lemma\rm}
\newtheorem{proposition} [theorem] {\sc  Proposition\rm}

\newtheorem{remark}[theorem]{\sc  Remark\rm}

\def\bproof{\noindent{\bf Proof.\;}}
\def\eproof{\hfill$\square$\medskip}

\newcounter{marnote}

\DeclareFontFamily{OT1}{rsfs}{}
\DeclareFontShape{OT1}{rsfs}{m}{n}{ <-7> rsfs5 <7-10> rsfs7 <10-> rsfs10}{}
\DeclareMathAlphabet{\mycal}{OT1}{rsfs}{m}{n}

\def\Id{{\rm Id}}

\def\tr{{\rm tr}}

\def\mcS{{\mycal{S}}}

\def\hh{{\rm \overline H}}

\def\Ss{{\mathbb S}}

\def\be{\begin{equation}}
\def\ee{\end{equation}}

\def\Id{{\rm Id}}

\def\tr{{\rm tr}}

\def\mcS{{\mycal{S}}}

\def\mcE{{\mycal E}}

\def \f {\varphi}

\newcommand{\defeq}{\stackrel{\scriptscriptstyle \text{def}}{=}}
\numberwithin{equation}{section}
\begin{document}

\title{Stability of the melting hedgehog in the Landau-de Gennes theory of nematic liquid crystals}
\author{Radu Ignat\thanks{Institut de Math\'ematiques de Toulouse, Universit\'e
Paul Sabatier, b\^at. 1R3, 118 Route de Narbonne, 31062
Toulouse, France. Email: Radu.Ignat@math.univ-toulouse.fr
}~, Luc Nguyen\thanks{Mathematical Insitute and St Edmund Hall, University of Oxford, Radcliffe Observatory Quarter, Woodstock Road, Oxford OX2 6GG, United Kingdom. Email: luc.nguyen@maths.ox.ac.uk}~, Valeriy Slastikov\thanks{School of Mathematics, University of Bristol, Bristol, BS8 1TW, United Kingdom. Email: Valeriy.Slastikov@bristol.ac.uk}~ and Arghir Zarnescu\thanks{University of Sussex, Department of Mathematics, Pevensey 2, Falmer, BN1 9QH, United Kingdom. Email: A.Zarnescu@sussex.ac.uk}}
\date{}

\maketitle
\begin{abstract}
We investigate stability properties of the radially symmetric solution corresponding to the vortex defect (so called ``melting hedgehog") in the framework of the Landau - de Gennes model of nematic liquid crystals.
We prove local stability of the melting hedgehog under arbitrary $Q$-tensor valued perturbations in the temperature regime near the critical {\it supercooling temperature}.  As a consequence of our method, we also  rediscover the  loss of stability of the vortex defect in the {\it deep nematic} regime.
\end{abstract}


\section{Introduction}
\label{sect:intropar}
In this paper we consider a variational Landau-de Gennes model for nematic liquid crystals. According to this model, the state of the liquid crystalline system is described 
at any point in space by a $Q$-tensor (i.e., a real three-by-three traceless symmetric matrix) and corresponds to stable critical points of Landau-de Gennes energy functional.  Our aim is to investigate the radially symmetric critical point of this functional that corresponds to the vortex defect, which is called {\it ``the melting hedgehog''} and has the following form:
\begin{equation}\label{def:meltinghedgehog}
H(x):=u(|x|)\left(\frac{x}{|x|}\otimes \frac{x}{|x|} -\frac{1}{3}Id\right), \quad x\in \RR^3\setminus\{0\},
\end{equation} where $u:\RR_+\to \RR$ with $u(0)=0$ is the radial scalar profile.

We will prove an important, both physically and mathematically, result about the local stability of $H(x)$ with respect to {\it arbitrary $Q$-tensor valued perturbations} in the regime where the temperature of the liquid crystal system is close to a critical temperature of isotropic-nematic phase transition, called the {\it ``supercooling temperature"}. Our approach consists in developing new tools for studying the second variation of the Landau-de Gennes functional, which might be of independent interest. We highlight that our stability result relies on the fine properties
of the scalar profile $u$ of the melting hedgehog that we proved in \cite{ODE_INSZ}. As a consequence of our approach, 
we also rediscover the local instability of $H(x)$ in the case of very low temperatures (called {\it ``the deep nematic regime"}), result first proved by Gartland and 
Mkaddem in \cite{mg}. 
 
In the next subsection we provide a physical motivation of this study
and the connection with other defect patterns. The more mathematically minded reader can skip directly to subsection~\ref{sec:math} where the notation and the mathematical results are presented. 

\subsection{Physical background and motivation}

Some of the most fascinating and striking achievements of modern material science are based on the use of materials with complex microstructure. In many cases the appearance of microstructure is associated with the presence of defects, i.e. localized regions where the material behavior is drastically different from the most common one. The complexity of defect patterns is intrinsically related to the structure of the material. One of the examples of technologically successful materials is provided by nematic liquid crystals -- the material of choice for constructing displays.  In liquid crystals the topological defects are easily observable experimentally and play an important role in the formation of the microstructure. There exist numerous  experimental and theoretical studies of defects in nematics (e.g., \cite{BallZar, chandra, cladiskle, dgdefects, klelavr, KraVirga, mg2, Saupe}). 
However there is still a lack in the proper understanding of the appearance and structure of defects in nematic liquid crystals.

Theoretical investigations of defect patterns are usually done within a framework of one of the current liquid crystal theories. 
We mention here only
the following macroscopic theories  attempting to describe spatially inhomogeneous liquid crystal systems: Oseen-Frank theory \cite{Frank}, Ericksen theory \cite{ericksen}, Onsager-type theory (see e.g. \cite{onsager-spatial_variation}), and Landau-de Gennes theory (see e.g. \cite{dg}). The state of the liquid crystal system is described by a macroscopic order parameter specific to each theory. The central object of study is a free energy functional depending on the order parameter, and equilibrium configurations of liquid crystal system correspond to local minimisers of the free energy. 

Let us describe briefly the order parameters used in the above theories. The Oseen-Frank theory describes the molecular alignment by objects with two degrees of freedom, namely unit-vectors $n(x)\in \Ss^2$ at any point $x$ in space. Ericksen theory uses, in addition to the unit-vector $n(x)$, a scalar parameter $s(x)\in [-\frac{1}{2}, 1]$, so in total three degrees of freedom. Onsager-type theory uses a vector field in $\RR^5$ that for critical points of the energy is in one-to-one correspondence with probability density of molecules orientations. Finally, Landau-de Gennes theory uses $Q$-tensors -- elements from the set of real three-by-three traceless symmetric matrices. We note that one can ``embed'' the Oseen-Frank and Ericksen theories in the more complex framework of the Onsager-type and  Landau-de Gennes theories. 

We are interested in studying {\it point defects} in the Landau-de Gennes theory. 
Although point defects are relevant for all above macroscopic theories, their definition varies in these theories.  In the Oseen-Frank theory, defects are defined as discontinuities of the order parameter $n$, where point defects with universal structure $n(x) = \frac{x}{|x|}$ are the only stable defects allowed in this theory (see \cite{BrezisCoronLieb, Hardt-Kinder-Lin}). In the Ericksen theory, point defects have the form $s(|x|) \frac{x}{|x|}$ so that the ``melting points'' correspond to the zeros of a scalar parameter $s$ (see \cite{Lin-Liu}). In Onsager-type theories, point defects correspond to $\delta$-functions (see \cite{Fat-Slas}). The situation is more complicated in Landau de-Gennes theory as there are different types of defects and several definitions of them; for instance, defects can be interpreted as points where two eigenvalues of the $Q$-tensor are equal (see \cite{KraVirga}), or where there is a discontinuity in the corresponding eigenframes (see \cite{dgdefects}). However, there is a natural extension of the Ericksen melting point defect to Landau-de Gennes theory. It has the structure of a ``melting hedgehog" $H(x)$ defined in \eqref{def:meltinghedgehog} with a defect located at the point $x=0$ where parameter $u$ vanishes. In this paper we investigate the stability of the melting hedgehog $H(x)$ in various temperature regimes within the framework of  Landau-de Gennes theory.
 
\subsection{Main results} 
\label{sec:math}
We consider the following Landau-de Gennes energy functional: 
  \begin{equation}
\mcF[Q; \Omega]:= \int_{\Omega} \Big[\frac{1}{2}|\nabla Q|^2 + f_{bulk}(Q)\Big]\,dx, \quad Q \in H^1(\Omega, \mcS_0), 
\label{L-DGFunc}
\end{equation} where $\Omega \subset \RR^3$ is a domain and the space of $Q$-tensors is denoted by $$\mcS_0:=\{ Q\in \RR^{3\times 3},\, Q=Q^t,\,\tr (Q)=0\}.$$ The potential $f_{bulk}$ accounts for the bulk effects and has the following form \footnote{ It is well-known that this form of a bulk energy density allows for multiple local minima and a first order nematic-isotropic phase transition \cite{dg,Virga}.}:
\begin{equation}
f_{bulk}(Q): = -\frac{a^2}{2}|Q|^2 - \frac{b^2}{3}\tr(Q^3) + \frac{c^2}{4}|Q|^4,
\label{BulkEDensity}
\end{equation}
where $a^2, b^2, c^2>0$ are material constants and $$|Q|^2:=\tr(Q^2).$$

A critical point $Q(x)$ of the functional $\mcF[\cdot; \Omega]$ satisfies the following Euler-Lagrange equation:
\begin{equation}
\Delta Q = -a^2\,Q - b^2[Q^2 - \frac{1}{3}|Q|^2\,Id] + c^2\,|Q|^2\,Q \quad \textrm{in } \Omega,
\label{L-DG::ELEqn}
\end{equation}
where the term $\frac{1}{3}b^2\,|Q|^2\,Id$ is a Lagrange multiplier that accounts for the tracelessness constraint. It is well known that solutions of \eqref{L-DG::ELEqn} are smooth (see for instance \cite{Ma-Za}).

\begin{remark}
\label{rem_scaling}
We should point out now that although the equation \eqref{L-DG::ELEqn} seems to depend on three parameters $a^2,b^2$ and $c^2$ there is only one independent parameter  in the problem, which can be chosen to be $a^2$, for fixed $b^2$ and $c^2$. This is because one can use a rescaling of $Q$ and its domain of definition $\Omega$ to fix two parameters $b^2$ and $c^2$.
Indeed, if $Q(x)$ is a solution of \eqref{L-DG::ELEqn} on the ball $\Omega:=B_R(0)$ centered at the origin, then the rescaled map $Q_{\lambda,\mu}(x):=\lambda Q(\frac{x}{\mu})$ solves
$$\Delta Q_{\lambda,\mu} = -\frac{a^2}{\mu^2}\,Q_{\lambda,\mu} - \frac{b^2}{\mu^2\lambda} [Q_{\lambda,\mu}^2 - \frac{1}{3}|Q_{\lambda,\mu}|^2\,Id] + \frac{c^2}{\mu^2\lambda^2}\,|Q_{\lambda,\mu}|^2\,Q_{\lambda,\mu},
$$
on the rescaled ball $B_{\frac{R}{\mu}}(0)$. Thus, given any $b_0^2$ and $c_0^2$ we can find a unique pair of  $\lambda>0$ and $\mu>0$ so that $\frac{b^2}{\mu^2\lambda}=b_0^2$ and $ \frac{c^2}{\mu^2\lambda^2}=c_0^2$, of course at the price of modifying $a^2$ and the radius of the domain. 

In what follows we will consider \emph{the coefficients $b^2$ and $c^2$ as fixed and vary the parameter $a^2$} that has a physical interpretation of a {\it reduced temperature} \footnote{Here we are making an assumption on the reduced temperature, which ensures that the isotropic phase $Q = 0$ is not a local minimizer of the bulk potential $f_{bulk}$. In contrast, if one replaces $a^2$ by $-a^2$ in $f_{bulk}$, the isotropic phase is locally stable.} (see \cite{mg2}). Alternatively, one could have fixed for example $a^2$ and $c^2$ and vary $b^2$, and we will see that indeed this makes sense from a mathematical point of view, when analyzing certain limit regimes. 
\end{remark}

We are interested in studying  {\it radially-symmetric} solution of \eqref{L-DG::ELEqn} in the whole space. We say that a map $Q\in H^1(\RR^3, \mcS_0)$ is radially-symmetric if
\begin{equation}
Q(\Rot x) = \Rot \,Q(x)\,\Rot^t \text{ for any } \Rot\in SO(3) \text{ and a.e. } x\in\RR^3.
\label{RadSymDef}
\end{equation}
In fact, such a map $Q(x)$ has only one degree of freedom (see \cite{ODE_INSZ}), namely,
there exists a continuous radial scalar profile $s: (0, +\infty) \rightarrow \RR$ such that 
$Q(x) = s(|x|)\,\hh (x)$
for a.e. $x \in \RR^3\;$,
where $\hh$ is called {\it hedgehog}
\be
\label{hedge}
\hh(x)=\frac{x}{|x|}\otimes \frac{x}{|x|} -\frac{1}{3}Id\ee and
the radial scalar profile $s$ of $Q$ is given by  
$s(|x|) = {3 \over 2} \tr (Q(x) \hh (x))$ for a.e. $x \in \RR^3$.

We will focus on the stability of the radially symmetric solution $H(x)$ of the system \eqref{L-DG::ELEqn}, i.e., the melting hedgehog defined in \eqref{def:meltinghedgehog} where the radial scalar profile $u:\RR_+\to \RR$ of $H(x)$ is a solution of the following boundary value problem:
\begin{align}
&u''(r) + \frac{2}{r}\,u'(r) - \frac{6}{r^2}\,u(r) = F(u(r)),
\label{RS::ODE0}\\
&u(0)=0, \quad \lim_{r\to\infty} \, u(r) =s_+. \label{BC}
\end{align} Here, $s_+\defeq\frac{b^2+\sqrt{b^4+24a^2c^2}}{4c^2}$ and \footnote{Here the value of $s_+$ is chosen so that the boundary condition at infinity for the melting hedgehog, namely $s_+\left(\frac{x}{|x|}\otimes \frac{x}{|x|}-\frac{1}{3}Id\right)$ is a global minimum of $f_{bulk}$ for every $x\neq 0$ (see \eqref{lim_manif}).} 
\begin{equation}\label{def:Fstandard}
F(u(r)):=
-a^2\,u(r) - \frac{b^2}{3}\,u(r)^2 + \frac{2c^2}{3}\,u(r)^3.
\end{equation} 

Our main result concerns the local stability of the melting hedgehog. We highlight that $H$ is a critical point of $\mcF(\cdot, \RR^3)$, but has infinite energy, i.e., $\mcF(H, \RR^3)=\infty$. Therefore, 
the second variation of the functional $\mcF$ at the point $H$ in the direction $V\in H^1(\RR^3;\mcS_0)$ is defined as \footnote{Equivalently, for every $V\in C_c^\infty(\RR^3;\mcS_0)$, one defines
$
\mcQ(V)
	= \frac{d^2}{dt^2}\Big|_{t = 0} \mcF[H + t\,V;\Omega]$
where the domain $\Omega$ is any bounded open set containing the support of $V$.}
\begin{align}
\nonumber
\mcQ(V)&= \frac 1 2\frac{d^2}{dt^2}\bigg|_{t=0} \int_{\RR^3} \Big[\frac{1}{2}|\nabla (H + t\,V)|^2 + f_{bulk}(H + t\,V)-\frac{1}{2}|\nabla H|^2 - f_{bulk}(H)\Big]\,dx\\
\label{secvar}
&= \int_{\RR^3} \Big[\frac{1}{2}|\nabla V|^2 + g(x,V)
\Big]\,dx 
\end{align}
where
\begin{equation}
g(x,V) := \big(-\frac{a^2}{2} + \frac{c^2u^2}{3}\big)|V|^2 - b^2\,u\,\tr(\bar H\,V^2) + c^2\,u^2\,\tr^2(\bar H\,V).
\label{Eq:gxV}
\end{equation}
The local stability of the melting hedgehog is formulated in the theorem below (already announced in our note \cite{INSZ_CRAS}):

\begin{theorem}
\label{thm:main} Let $b^2$ and $c^2$ be fixed positive constants. Then there exists $a^2_0>0$ (depending on $b^2$ and $c^2$)  such that for all $a^2<a^2_0$ the radially symmetric solution $H$ is locally stable in $H^1(\RR^3;\mcS_0)$, meaning that $\mcQ(V)\geq 0$ for all $V \in H^1(\RR^3;\mcS_0)$. Moreover 
$\mcQ(V)=0$ if and only if $V \in \{ \partial_{x_i} H \}_{i=1}^3$, i.e. the kernel of the second variation coincides with translations of $H(x)$. 

Also, there exists $a_1^2>0$ (depending on $b^2$ and $c^2$) so that for any $a^2>a_1^2$ there exists $V_*\in C_c^\infty(\RR^3;\mcS_0)$ such that $\mcQ(V_*)<0$. Moreover, any such $V_*$ cannot be purely uniaxial (i.e., $V_*(x)$ has three different eigenvalues for some point $x\in \RR^3$).
\end{theorem}

\begin{remark}
In Theorem \ref{thm:main}, when $a^2 > a_1^2$, $V_*$ can be chosen axially symmetric, i.e. there exists an axis $L$ through the origin such that if $\mathcal{R}$ is any rotation matrix about $L$, there holds
\[
V_*(\mathcal{R}x) = \mathcal{R}\,V_*(x)\,\mathcal{R}^t \text{ for all } x \in \RR^3.
\]
In particular, the melting hedgehog is unstable within the class of axially symmetric $Q$-tensors. Therefore, $\mcF$ has an axially symmetric critical point which does not coincide with the melting hedgehog.
\end{remark}

As a consequence, we obtain the following result concerning the local minimality of the melting hedgehog.
\begin{corollary}\label{Cor:MHLocMin} Let $b^2$ and $c^2$ be fixed positive constants. Then there exists $a^2_0>0$ (depending on $b^2$ and $c^2$) such that
for $a^2 < a^2_0$, the melting hedgehog $H$ is locally minimizing, i.e. for any open bounded subset $\Omega$ of $\RR^3$, there exists a constant $\epsilon_0 = \epsilon_0(\Omega) > 0$ such that $\mcF[H + V;\Omega] \geq \mcF[H; \Omega]$ for all $V \in H^1_0(\Omega;\mcS_0)$ with $\|V\|_{H^1} \leq \epsilon_0$.
\end{corollary}

Another consequence of our analysis of the second variation is the following result:
\begin{proposition}\label{Cor:RestrictiveStability}
For all $a^2 \geq 0$, the melting hedgehog $H$ is locally stable with respect to uniaxial perturbations, i.e., $\mcQ(V)\geq 0$ for every $V\in H^1(\RR^3;\mcS_0)$ such that $V(x)$ has two equal eigenvalues at each point $x\in \RR^3$.
\end{proposition}

 The stability of the  melting hedgehog has already attracted substantial interest in the physical and mathematical literature on liquid crystals. We mention here the early studies of R. Rosso and E.G. Virga \cite{Rosso-Virga} who obtained some local stability results under restricted classes of perturbations. Later S. Mkaddem and E.C. Gartland, \cite{mg,mg2} investigated numerically and analytically the stability properties of the melting hedgehog \eqref{def:meltinghedgehog} considering only axially symmetric perturbations $V$. Using a combination of numerical and analytical arguments, they obtained the local instability of the melting hedgehog for large enough values of the parameter $a^2$ (by providing a suitable perturbation $V$ that makes $\mcQ(V)< 0$). Their numerical studies suggest that for small enough $a^2$ the melting hedgehog \eqref{def:meltinghedgehog} is locally stable, a result that we prove rigorously in Theorem~\ref{thm:main}.

\subsection{Comparison with Ginzburg-Landau systems}
Let us compare the Landau-de Gennes (LdG) model with the $2D$ and $3D$ Ginzburg-Landau (GL) systems \cite{BBH, Mironescu-radial, Mil-Pis}. Both LdG and GL energies consist of two terms: a Dirichlet energy and a potential energy. The $N$-dimensional GL energy is defined for $N$-dimensional vector fields (here, $N=2,3$) and the minima of the well potential is a co-dimension one manifold, e.g., the unit sphere $\Ss^{N-1}\subset \RR^N$. 
However, our $3D$ LdG energy is defined for maps with values into the five-dimensional space $\mcS_0$ of $Q$-tensors and the global minimizing $Q$-tensors of the bulk potential $f_{bulk}$, as defined in \eqref{BulkEDensity}, form a $2D$ manifold ${\cal M}\subset \mcS_0$ (the limit target manifold), homeomorphic to the projective plane $\RR P^2$ 
\be
\label{lim_manif}
{\cal M}=\Big\{Q\in \mcS_0\, :\, Q=s_+\bigg(n\otimes n -\frac{1}{3}Id\bigg), \, n\in \Ss^2\Big\}.
\ee 
Thus, one might think at first that similar stability behavior holds for the hedgehog for the LdG and the vortex defect for the $3D$ GL (as maps defined from $\RR^3$ into a $2D$ manifold), see \cite{Gustafson}. Nevertheless, there is still a significant difference to $3D$ GL since variations in LdG are allowed in $5D$ and the target limit manifolds ($\RR P^2$ for LdG and $\Ss^2$ for $3D$ GL) have different topologies. 
We highlight that the situation would be even more different for the LdG model in the limit regime $b^2 = 0$: in that case,
the global minimizing $Q$-tensors of $f_{bulk}$ form a $4D$ sphere $\{Q\in \mcS_0\,:\, \tr(Q^2)=a^2/c^2\}$. In fact, it was proved in \cite{mg} that the melting hedgehog is unstable if $b^2$ is small (corresponding to $a^2$ large after the rescaling, see Remark \ref{rem_scaling}). 
Our major contribution in this paper is to go beyond the Ginzburg-Landau type of methods and develop a systematic approach that is capable of dealing not only with the regime $a^2 \to \infty$ but also with the very different and much more challenging case when $a^2$ is small.

\subsection{Sketch of the method}
An essential step in our analysis relies on the fine properties of the scalar profile $u$. Its existence and uniqueness, together with a refined study of its qualitative properties was accomplished in \cite{ODE_INSZ}. The main results needed in this paper are summarised in the following: 
  
\begin{theorem}[\cite{ODE_INSZ}]\label{thm:ODE}
There exists a unique solution $u\in C^\infty([0,\infty))$ of \eqref{RS::ODE0} with the right-hand side $F(u)$ as in \eqref{def:Fstandard} and boundary conditions \eqref{BC}. This solution is positive and strictly increasing in $(0,\infty)$.

Furthermore, if we denote $w(r):=\frac{ru'(r)}{u(r)}$ and $f(u)=\frac{F(u)}{u}$ then for every $r\in (0, \infty)$ we have:
\be
u(r) > u_0(r), \hbox{ where $u_0$ solves \eqref{RS::ODE0}, \eqref{BC} for $a^2=0$}, 
\ee
\be
u'(0)=0 \hbox{ and } u''(0) \geq u''_0(0)>0,
\label{Eq1.14}
\ee
\be \label{rel:uprimes}
u'' + (-\frac{3u'}{u} + \frac{5}{r})u' \geq 0,
\ee
\be \label{rel:relfs}
2a^2 + \frac{b^2}{3} u > -\frac{2}{w}\,f(u),
\ee
\be
\label{rel:fws}
\frac{3}{r^2}(w-2)(w+ 1)< f(u)
< \frac{1}{r^2}(w-2)(2w + 3),
\ee
\begin{equation}
\label{rel:u'u}
0<w(r) < 2,
\end{equation} 
\begin{equation}
\label{rel:u'u1}
\frac{1}{w(r)} > \frac{u''(0)}{4s_+} r^2.
\end{equation}
\end{theorem}

We now provide a brief overview of the  arguments for proving our stability results that are expanded and detailed in Sections \ref{sec:prelim}, 
\ref{sec:SepVar} and \ref{sec:stability}.
We split investigation of the second variation $\mcQ(V)$ into several steps. The first idea is to represent any perturbation $V$ as a linear combination with respect to certain {\it pointwise orthogonal} frame in the set $\mcS_0$ of traceless symmetric $3 \times 3$ matrices. Using this approach we obtain the representation of the second variation $\mcQ(V)$ in terms of five scalar functions $w_i$, $0 \leq i \leq 4$ (see \eqref{energy}).

It is clear that each of the five functions $\{w_i\}_{0 \leq i \leq 4}$ depend on the three spherical coordinates characterizing a point $x\in \RR^3$: $r\in (0,\infty),\ \theta\in [0,\pi], \ \varphi\in [0,2\pi)$. Our approach consists in ``removing" the dependence on ``spherical variables" $(\theta, \varphi)$ by using at each step a special basis decomposition of functions $w_i$.

More precisely, we use the standard Fourier decomposition in variable $\varphi$ of $w_i$, $0\leq i \leq 4$ and as a consequence, represent the second variation as a series $\mcQ(V)=\sum_k \mcQ(V_k)$, where $\{V_k\}_{k\geq 0}$ are the projection of $V$ onto the corresponding basis induced subspaces of $\mcS_0$ (see \eqref{Vk}). Then we show that all the 
modes $\mcQ(V_k)\geq 0$ are nonnegative for $k>2$ provided $\mcQ(V_k)\geq 0$ for $k=0, 1, 2$. Therefore we just have to investigate the sign of the first three modes $\mcQ(V_0)$, $\mcQ(V_1)$, $\mcQ(V_2)$. Furthermore, we show in Proposition \ref{Phi012} that due to some symmetries, our approach reduces to the study of the sign of the \emph{three functionals} $\Phi_k$, $k=0,1,2$ depending only on \emph{three $\varphi$-independent variables} $v_0$, $v_2$, $v_4$ (see \eqref{Phidef}). However the problem is still far from being simple since the $v_{2k}$'s depend on $r$ and $\theta$ variables and the sign of $\Phi_k$ is not transparent due to coupling terms $v_i v_j$ and $v_i \partial_\theta v_j$. 

The next step consists in eliminating the $\theta$-dependence and simplifying the problem even further.
The analysis here is quite delicate and involves a decomposition of $\{v_{2k}(r, \theta)\}_{k=0,1,2}$ using \emph{suitably adapted bases}. We first use the formal asymptotics when $r\to\infty$ to gain the intuition  and discover these bases as solutions to specific eigenvalue problems for suitable elliptic operators with non-constant coefficients (see Section~\ref{sec-heuristics}). Expanding $v_{2k}$, $0\leq k \leq 2$  in these systems of eigenvectors, we decompose $\Phi_k = \pi \sum_{j\geq 0} \Phi_{k,j}$, $0 \leq k \leq 2$
where each functional $\Phi_{k,j}$ is defined on three functions depending {\it only on $r$}. The last identity, i.e. an orthogonality relation \emph{for a system}, is \emph{unexpected} and dramatically simplifies our study. 

 At this point, it is enough to study only the sign of $\Phi_{0,j}$ since $\{ \Phi_{0,j} \}_{j=0}^\infty \supset \{ \Phi_{k,j} \}_{j=0}^\infty$ for $k=1,2$. 
It is not so difficult to show that $\Phi_{0,j} \geq 0$ for $j \neq 2,3$. We have thus reduced the original problem of finding the sign of the second variation $\mcQ (V)$ to investigation of the sign of the functionals $\Phi_{0,j} (\tilde{w}_0, \tilde{w}_2, \tilde{w}_4)$, $j = 2,3$ defined in \eqref{new_form_phi} where $\tilde{w}_0, \tilde{w}_2, \tilde{w}_4$ depend only on $r$. This is a nontrivial problem due to presence of nonlinearities $f(u)$, $\hat f(u)$, $\tilde f(u)$ (defined at \eqref{fu}, \eqref{fhat}, \eqref{ftil}) and coupling terms in the expression of $\Phi_{0,j}$. We use suitable Hardy-type decompositions of functions $\tilde{w}_0(r)$, $\tilde{w}_2(r)$, $\tilde{w}_4(r)$ to eliminate the above nonlinearities  together with the fine properties of the profile $u$ in Theorem \ref{thm:ODE} in order to conclude that $\Phi_{0,j}\geq0$ ($j = 2,3$) when the parameter $a^2$ is small enough.

\section{A special basis of $Q$-tensors. An instability result.}
\label{sec:prelim}

Let $H$ be the melting hedgehog \eqref{def:meltinghedgehog} where the radial scalar profile $u$ is the solution of \eqref{RS::ODE0} \& \eqref{BC} given by Theorem \ref{thm:ODE}. Consider the second variation $\mcQ(V)$ of $\mcF(\cdot, \RR^3)$ at $H$ in the direction $V\in C^\infty_c(\RR^3;\mcS_0)$, which is the quadratic functional given in \eqref{secvar}. Note that since $u$ is bounded, then there exists a constant $C>0$ (depending on the parameters $a^2, b^2, c^2>0$) such that $|g(x,Q)| \leq C\,|Q|^2$ for every $Q\in \mcS_0$ and $x\neq 0$ and so $\mcQ(\cdot)$ extends to a bounded quadratic functional on $H^1(\RR^3;\mcS_0)$.

As mentioned in the introduction, the goal of the paper is to investigate the sign of the second variation $\mcQ(\cdot)$, especially in the following two regimes: $a^2$ large and $a^2$ small, respectively. We start by introducing in subsection \ref{basisdec} a locally adapted frame in the space $\mcS_0$ of $Q$-tensors near the melting hedgehog $H$ and the expression of the second variation in such frame. We then follow up in subsection \ref{instab} with some simple 
analysis that leads to the proof of instability of $H$ stated at Theorem \ref{thm:main} when $a^2$ is sufficiently large.

\subsection{Decomposition in a special basis of $Q$-tensors} \label{basisdec}

The idea is to represent a perturbation $V\in C^\infty_c(\RR^3, \mcS_0)$ as a linear combination with respect to a certain orthogonal frame in the set $\mcS_0$ of traceless symmetric $3 \times 3$ matrices. In the sequel, it is convenient to use spherical coordinates for $x\in \RR^3\setminus\{0\}$:
$$
x= r  ( \sin\theta\,\cos\varphi , \sin\theta\,\sin\varphi , \cos\theta ) , \quad r>0, \, \theta\in [0, \pi], \, \varphi \in [0, 2\pi).
$$
We define the following orthonormal basis in $\RR^3$:
$$
\begin{array}{lll}
 n&=& ( \sin\theta\,\cos\varphi , \sin\theta\,\sin\varphi , \cos\theta )=\frac{x}{|x|} , \\
m &=& (\cos\theta\,\cos\varphi , \cos\theta\,\sin\varphi , -\sin\theta ), \\
p &=& (\sin\varphi , -\cos\varphi , 0 ).
\end{array}
$$
Using this basis, we can also define an orthogonal frame in $\mcS_0$ as
\begin{align}
\label{def:Ei}
E_0=\bar H=  n \otimes n& - {1 \over 3} \IdMa, \ E_1 = n \otimes p + p \otimes n, \ E_2 = n \otimes m + m \otimes n,\nonumber\\
&E_3 = m \otimes p + p \otimes m, \ E_4 = m \otimes m - p \otimes p.
\end{align}
Note that $E_i=E_i(\theta, \varphi)$ are independent of $r$; moreover, $\{E_1, E_2\}$ (resp., $\{E_0, E_3, E_4\}$) form a basis of the tangent space (resp., the normal space) of the limit target manifold ${\cal M}$ at $s_+\bar{H}$ (see definitions \eqref{hedge} and \eqref{lim_manif}), see \cite{Ngu-Zar}. It is clear that any $\mcS_0$-valued map $V$ can be represented as the following linear combination:
\be \label{defV}
V(x) = \sum_{i=0}^4 w_i (r,\theta, \varphi) E_i(\theta, \varphi), \quad x\neq 0,
\ee for some scalar functions $w_0, \ldots, w_4$.

We want now to decompose the second variation $\mcQ(V)$ in this frame. First, note that $E_i \cdot E_j =\tr (E_i E_j^t)=0$ for $i \neq j$ and
\[
|E_0|^2 = \frac{2}{3}, \ |E_i|^2 = 2 \ \hbox{ for } i=1,\dots,4.
\]
A straightforward calculation gives for every $x\neq 0$:
\begin{align*}
&|V|^2 = {2 \over 3} w_0^2 + 2 \sum_{i=1}^4 w_i^2, \qquad \tr^2 (\hh V) = {4 \over 9} w_0^2, \\
& \tr(\hh V^2) = {2 \over 9} w_0^2 + {1 \over 3} w_1^2 + {1 \over 3} w_2^2 -{2 \over 3} w_3^2 - {2 \over 3} w_4^2 .
\end{align*}
Recalling definition \eqref{Eq:gxV}, 
we thus obtain:
$$
g(x,V(x)) = {1 \over 3} w_0^2 \hat f (u(r)) + (w_1^2 + w_2^2) f(u(r)) + (w_3^2 + w_4^2) \tilde f(u(r)), \quad x\neq 0,
$$
where 
\begin{align}
& f(u(r)) = \frac{F(u(r))}{u(r)} = - a^2 - {b^2 u(r) \over 3} +\frac{2 c^2 u(r)^2}{3}
	,\label{fu} \\
&\hat f(u(r)) = F ' (u(r)) = - a^2 - {2 b^2 u(r) \over 3} +2 c^2 u(r)^2
	, \label{fhat}\\
&\tilde f(u(r)) = -a^2 + \frac{2b^2u(r)}{3} + \frac{2c^2u(r)^2}{3}
	.\label{ftil}
\end{align}

In order to see how the kinetic term $|\nabla V|^2$ breaks up with respect to the basis $\{E_i\}_{0\leq i \leq 4}$, we first compute the derivatives of the basis elements. We note that, at $x\neq 0$, there hold: $\partial_\theta n = m, \partial_\theta m = -n, \partial_\theta p = 0,
\partial_\varphi n = -\sin\theta\,p, \partial_\varphi m = -\cos\theta\,p,  \partial_\varphi p = \sin\theta\,n + \cos\theta \, m$. As a consequence, we find that
\begin{align}
\nonumber
&\partial_\theta E_0 = E_2, \ \partial_\theta E_1 = E_3,\ \partial_\theta  E_2 = - 3E_0 + E_4,\ \partial_\theta E_3 = -E_1,\ \partial_\theta E_4 = - E_2,\\
\nonumber
&\partial_\varphi E_0 = - \sin\theta\,E_1,\ \partial_\varphi E_1= 3\sin\theta\,E_0 + \cos\theta\, E_2 + \sin\theta\,E_4,\\
\nonumber
&\partial_\varphi E_2 = -\cos\theta E_1 - \sin\theta E_3,\ \partial_\varphi E_3 = \sin\theta\,E_2 + 2\cos\theta\,E_4,\ \partial_\varphi E_4 = -\sin\theta E_1 - 2\cos\theta E_3.
\end{align}
This implies
\begin{align*}
\partial_\theta V
	&= (\partial_\theta w_0 - 3  w_2)E_0
		+ (\partial_\theta w_1 - w_3)E_1
		+ (\partial_\theta  w_2 + w_0 - w_4) E_2\\
	&\quad + (\partial_\theta w_3 + w_1)E_3
		+ (\partial_\theta w_4 +  w_2)E_4
	,\\
\partial_\varphi V
	&= (\partial_\varphi w_0 + 3\sin\theta\,w_1)E_0
		+ (\partial_\varphi w_1 -\sin\theta\, w_0 - \cos\theta\,  w_2 - \sin\theta\,w_4)E_1\\
	&\quad + (\partial_\varphi  w_2 + \cos\theta\,w_1 + \sin\theta\,w_3) E_2
		+ (\partial_\varphi w_3 - \sin\theta\,  w_2 - 2\cos\theta\, w_4)E_3\\
	&\quad + (\partial_\varphi w_4 + \sin\theta\, w_1 + 2\cos\theta\, w_3)E_4.
\end{align*}
From this as well as the orthogonality property of the $E_i$'s, we see that
\begin{align}
\mcQ(V)
	= &\int_0^\infty\,\int_0^{2\pi}\,\int_0^\pi \bigg\{\frac{1}{3}|\partial_r w_0|^2 + |\partial_r w_1|^2 + |\partial_r  w_2|^2 + |\partial_r w_3|^2 + |\partial_r w_4|^2\nonumber\\
		&+ \frac{1}{r^2}\Big[\frac{1}{3}(\partial_\theta w_0 - 3  w_2)^2
			+ (\partial_\theta w_1 - w_3)^2
			+ (\partial_\theta w_2 + w_0 - w_4)^2 \nonumber \\
			&\quad + (\partial_\theta w_3 + w_1)^2
			+ (\partial_\theta w_4 +  w_2)^2\Big] \nonumber \\
		&+ \frac{1}{r^2\sin^2\theta}\Big[\frac{1}{3}(\partial_\varphi w_0 + 3\sin\theta\,w_1)^2 + (\partial_\varphi w_1 -\sin\theta\, w_0 - \cos\theta\,  w_2 - \sin\theta\,w_4)^2 \nonumber \\
			&\quad+ (\partial_\varphi  w_2 + \cos\theta\,w_1 + \sin\theta\,w_3)^2
			+ (\partial_\varphi w_3 - \sin\theta\,  w_2 - 2\cos\theta\, w_4)^2 \nonumber \\
			&\quad + (\partial_\varphi w_4 + \sin\theta\, w_1 + 2\cos\theta\, w_3)^2\Big] \nonumber \\
		&+ \frac{1}{3}\,\hat f(u)\,w_0^2 + f(u) \left(  w_1^2 + w_2^2 \right) + \tilde f(u) \left( w_3^2 + w_4^2 \right) \bigg\}\,r^2\,\sin\theta\,d\theta\,d\varphi\,dr
	.\label{energy}
\end{align}

\begin{remark}We already know that $\mcQ(V)$ is finite for every $V\in H^1(\RR^3;\mcS_0)$. Due to expression \eqref{energy}, it is important to
investigate under which condition on $\{w_i\}_{0 \leq i \leq 4}$, the associated $V$ in \eqref{defV} belongs to $H^1$. In fact, we note that for given $w_0, \ldots, w_4 \in L^2(\RR^3)$, setting $V= \sum_{i=0}^4 w_i\,E_i$, then $V\in H^1(\RR^3;\mcS_0)$ if and only if
\begin{multline}
|\nabla w_0| + |\nabla w_1 - \frac{\cos\theta}{r^2\sin^2\theta} w_2 \partial_\varphi| + |\nabla w_2 + \frac{\cos\theta}{r^2\sin^2\theta} w_1 \partial_\varphi|\\
	 + |\nabla w_3 - \frac{2\cos\theta}{r^2\sin^2\theta} w_4 \partial_\varphi| + |\nabla w_4 + \frac{2\cos\theta}{r^2\sin^2\theta} w_3 \partial_\varphi| \in L^2(\RR^3).
	 \label{Eq:wiSpace}
\end{multline}
Here $\partial_\varphi$ denotes the vector $(-x_2, x_1, 0)$. 
In particular, \eqref{Eq:wiSpace} is fulfilled provided that
\begin{multline}
w_0 \in H^1(\RR^3) \text{ and } \\
	w_j \in H^1_L(\RR^3) := \Big\{w \in L^2(\RR^3): |\nabla w| + \frac{|w|}{r\sin\theta} \in L^2(\RR^3)\Big\}, j = 1, \ldots, 4.
	\label{Eq:wiSpaceX}
\end{multline}
However, \eqref{Eq:wiSpace} and \eqref{Eq:wiSpaceX} are not equivalent. For example, if $\chi \in C_c^\infty(B_2(0))$ is a smooth non-negative cut-off function with $\chi = 1$ in $B_1(0)$, then $w_0 = w_3 = w_4 = 0$, $w_1 = r \chi\cos\theta\cos\varphi $ and $w_2 = -r\chi\cos^2\theta\sin\varphi$ satisfy \eqref{Eq:wiSpace} but not \eqref{Eq:wiSpaceX}.
\end{remark}

\begin{remark}
In the analysis of $\mcQ(V)$, we will frequently use the sharp Hardy inequality in three dimensions without explicitly mentioning:
\[
\int_{\RR^3} |\nabla \psi|^2\,dx \geq \frac{1}{4}\int_{\RR^3} \frac{|\psi|^2}{|x|^2}\,dx \text{ for all } \psi \in H^1(\RR^3).
\]
\end{remark}


\subsection{Instability of the melting hedgehog in the regime $a^2$ large} \label{instab}

We analyze now the expression \eqref{energy} and show that there is a range of parameters $a^2$, $b^2$, and $c^2$ such that $\mcQ(V)$ becomes negative for certain $\{w_i\}_{0 \leq i \leq 4}$ in \eqref{defV}. Instead of jumping immediately to the expression of the test functions $\{w_i\}_{0 \leq i \leq 4}$, we would like to show the logic of the construction that will be helpful later. Therefore we start with the natural reduction of the problem to the {\it axially symmetric case}, i.e., 
$$\textrm{$w_i=w_i(r, \theta)$ are $\varphi$-independent.}$$ Moreover, we assume the simplest possible ansatz when {\it only one} of the functions $w_i$ is not identically equal to zero, i.e., we focus on the subspaces
$$\mcV_i:=\{w_i E_i\, :\, w_i=w_i(r, \theta)\}, \quad i=0, \dots, 4$$ 
(note that here and in the rest of the section we do not assume summation over repeated indices). Under these assumptions, the second variation becomes
\begin{align*}
\mcQ(V)
	=   a_i \mcQ_i (w_i) \textrm{ for } V=w_i(r, \theta) E_i(\theta, \f)\in \mcV_i,
\end{align*}
where $a_0 ={1 \over 3}$, $a_i =1$ for $0 \leq i \leq 4$, and the functionals $\mcQ_i$ are given by
\begin{align*}
\mcQ_0(w_0)
	&= 2 \pi\int_0^\infty \int_0^\pi \Big\{|\partial_r w_0|^2 + \frac{1}{r^2}|\partial_\theta w_0|^2 + \frac{6}{r^2}\,w_0^2 + \hat f(u) \,w_0^2\Big\}\,r^2\,\sin\theta\,dr\,d\theta
	,\\
\mcQ_1(w_1)
	&= 2 \pi\int_0^\infty \int_0^\pi \Big\{|\partial_r w_1|^2 + \frac{1}{r^2}|\partial_\theta w_1|^2 + \frac{1}{r^2}(4 + \frac{1}{\sin^2\theta})\,w_1^2 + f(u)\,w_1^2\Big\}\,r^2\,\sin\theta\,dr\,d\theta
	,\\
\mcQ_3(w_3)
	&= 2 \pi\int_0^\infty \int_0^\pi \Big\{|\partial_r w_3|^2 + \frac{1}{r^2}|\partial_\theta w_3|^2 + \frac{1}{r^2}(-2 + \frac{4}{\sin^2\theta})\,w_3^2 + \tilde f(u)\,w_3^2\Big\}\,r^2\,\sin\theta\,dr\,d\theta,
\end{align*}
and $\mcQ_2(w) = \mcQ_1(w)$, $\mcQ_4(w) = \mcQ_3(w)$ for $w=w(r, \theta)$. 

\begin{remark} 
It is readily seen that $\mcQ_0(w)$ is finite for $w(r, \theta)\in H^1(\RR^3)$ while $\mcQ_1(w), \ldots, \mcQ_4(w)$ are finite for $w(r, \theta) \in H^1_L(\RR^3)$ (see \eqref{Eq:wiSpaceX}); therefore, we extend $\{\mcQ_i\}_{1 \leq i \leq 4}$ for $w(r, \theta) \in H^1(\RR^3)$ by letting them take the value $+\infty$ if $\frac{w}{r\sin\theta} \notin L^2(\RR^3)$.
\end{remark}

If either one of the above functionals becomes negative for some test function, then the instability of the radial hedgehog is proved. We have the following result:
\begin{theorem} \label{thm:scalar}
For any $w \in H^1 (\RR^3)$  such that $w=w(r, \theta)$ is independent of $\varphi$-variable we have
\begin{enumerate}
\item (Stability over $\mcV_0$) $\mcQ_0(w) \geq 0$; 
\item (Stability over $\mcV_1$ and $\mcV_2$) $\mcQ_i(w) \geq 0$ for $i=1, 2$;
\item (Stability over $\mcV_3$ and $\mcV_4$ for small $a^2$) For fixed positive constants $b^2$ and $c^2$, there exists $a^2_0>0$ (depending on $b^2$ and $c^2$) 
such that if $a^2< a^2_0$ then $ \mcQ_i(w) \geq 0$ for $i=3, 4$;
\item  (Instability for large $a^2$) For fixed positive constants $b^2$ and $c^2$, there exists $a^2_1 >0$ (depending on $b^2$ and $c^2$) such that if $a^2 > a^2_1$ then we can construct $w_3 \in C_c^\infty(\RR^3)$ independent of the $\varphi$-variable with $\mcQ_3(w_3) <0$. 
\end{enumerate}
\end{theorem}

\begin{proof}
Without loss of generality, we can assume that $w = w(r,\theta) \in C_{c}^\infty(\RR^3 \setminus \{0\})$ when proving {\it 1., 2.} and {\it 3.} 


\noindent {\it Proof of 1.} We start by showing that  $\mcQ_0$  is positive definite, using an explicit computation to reveal a nice structure of the functional $\mcQ_0$ that will be used later. Let $w(r, \theta)\in C_c^\infty(\RR^3 \setminus \{0\})$. It's clear that
$$
\frac{\mcQ_0 (w)}{2\pi} \geq  \int_0^\pi \left( \int_0^\infty  \Big\{|\partial_r w|^2  + \frac{6}{r^2}\,w^2 + \hat f(u) \,w^2\Big\}\,r^2 \, dr \right)\,\sin\theta\,d\theta .
$$ 
We can represent $w(r,\theta)$ as $ w(r,\theta) = u'(r) \ringw (r, \theta)$, where $\ringw  \in C_{c}^\infty (\RR^3  \setminus \{0\}) $.
Since $u$ is a solution of the ODE \eqref{RS::ODE0}, we deduce that:
\[
u''' + \frac{2}{r}\,u'' - \frac{8}{r^2}\,u' + \frac{12}{r^3}u = F'(u)u'=\hat f(u)\,u'.
\]
Therefore, it is straightforward to obtain \footnote{This is an application of the more general Hardy-type decomposition described in Lemma~\ref{thm:hardy} (see Appendix).}
\begin{align*}
\frac{\mcQ_0 (w)}{2\pi}
	& \geq \int_0^\pi \left( \int_0^\infty \Big[|u''\ringw + u'\,\partial_r\ringw|^2 + \frac{6}{r^2}\,|u'|^2\ringw^2 + \hat f(u)\,|u'|^2\,\ringw^2\Big]\,r^2\,dr \right) \sin\theta \, d\,  \theta \\
	&= \int_0^\pi \left(  \int_0^\infty \Big[|u'|^2\,|\partial_r\ringw|^2 + ( - \frac{2}{r^2}\,u' + \frac{12}{r^3}u)u'\,\ringw^2\Big]\,r^2\,dr  \right) \sin\theta \, d\, \theta,
\end{align*}
where we used integration by parts in $r$-variable for the term $u'''\cdot(u'\,\ringw^2r^2)$. 
Using that $0<u'\leq 2u/r$ for $r>0$ (see \eqref{rel:u'u} in Theorem~\ref{thm:ODE}), we obtain 
\be\label{mcE0positive}
\mcQ_0 (w) \geq 4 \int_{\RR^3} \frac{w^2}{|x|^2} \, dx \geq 0, 
\ee 
which proves that $\mcQ_0$ is indeed positive definite over $\{w\in H^1(\RR^3)\, :\, w=w(r, \theta)\}$.

\medskip

\noindent {\it Proof of 2.} 
Now we proceed to investigate $\mcQ_1$ and show that it is also positive definite. We have the following inequality which is a consequence of Lemma~\ref{thm:hardy} (see Appendix): for any real $k > 0$ and $v\in H^1\big((0,\pi);\sin\theta d\theta\big)$,
\be \label{PWH}
\int_0^\pi \Big[|v'|^2 + \frac{k^2}{\sin^2\theta}\,v^2\Big]\,\sin\theta\,d\theta \geq (k^2 +k) \int_0^\pi v^2\,\sin\theta\,d\theta, 
\ee
and equality is achieved if and only if $v(\theta) = C \sin^k (\theta)$ for some constant $C$. Let now $w(r, \theta)\in C_{c}^\infty (\RR^3  \setminus \{0\})$.
Using \eqref{PWH} (for $k=1$), we obtain:
$$
\mcQ_1(w) \geq 2 \pi \int_0^\pi \left( \int_0^\infty \Big[ | \partial_r w|^2 + \big(\frac{6}{r^2} + f(u)\big)| w|^2\Big]\,r^2\,dr \right) \sin \theta \, d \theta .
$$
Using representation $w (r, \theta) =u (r) \ringw (r, \theta)$ with $\ringw \in C_{c}^\infty (\RR^3  \setminus \{0\})$, it follows by \eqref{RS::ODE0}: 
\begin{align*}
\mcQ_1(w)
	&\geq 2 \pi \int_0^\pi \left( \int_0^\infty \Big[|u' \,\ringw + u\,\partial_r \ringw|^2 + \frac{6}{r^2}\,u^2\,\ringw^2 + f(u)\,u^2\,\ringw^2\Big]\,r^2\,dr  \right) \sin \theta \, d \theta\\
	&= 2 \pi \int_0^\pi  \int_0^\infty u^2\,|\partial_r \ringw|^2\,r^2\,dr \,  \sin \theta \, d \theta, 
\end{align*}
where we used integration by parts in $r$-variable for the term $u''\cdot(u\,\ringw^2r^2)$.
Therefore we obtain the second assertion that $\mcQ_1(w) >0$ unless $w \equiv 0$.

\medskip

\noindent {\it Proof of 3.} 
We now prove that $\mcQ_3 (w)$ is positive definite on $\mcV_3$. Let $w(r, \theta)\in C_{c}^\infty (\RR^3  \setminus \{0\})$. Using \eqref{PWH} (for $k=2$) we obtain
\be \label{E3}
\mcQ_3(w) \geq 2 \pi \int_0^\pi \left( \int_0^\infty \Big[ | \partial_r w|^2 + \big(\frac{4}{r^2} + \tilde f(u)\big)| w|^2\Big]\,r^2\,dr \right) \sin \theta \, d \theta .
\ee
Using Theorem~\ref{thm:ODE} we know that $u \geq u_0$ on $[0, \infty)$ and therfore we obtain
\[
\frac{4}{r^2} + \tilde f(u) \geq \Big( \frac{4}{r^2} + \frac{2b^2}{3} u_0 + \frac{2c^2}{3} u_0^2\Big) - a^2\geq a^2_0-a^2,
\]
where $a^2_0>0$ is a minimum of the function  $\frac{4}{r^2} + \frac{2b^2}{3} u_0 + \frac{2c^2}{3} u_0^2$ on the interval $[0,\infty)$ and depends only on $b^2, c^2$ . Choosing $a^2 < a_0^2$ we conclude the proof.
\medskip

\noindent {\it Proof of 4.} 
Let us now show that instability appears for $a^2$ large. Namely, we need to construct $w_3(r,\theta) \in C_c^\infty(\RR^3)$ such that $\mcQ_3(w_3) < 0$. Recall that the equality in \eqref{PWH} (for $k=2$) is achieved for $\psi(\theta)=\sin^2 \theta$. Therefore, we restrict to $$w_3 (r, \theta) = w_*(r) \sin^2 \theta, \qquad w_* \in C_c^\infty(0,\infty).$$ By the previous computation in \eqref{E3}, we obtain that 
$$
\mcQ_3(w_3) = \frac{32 \pi}{15} 
\int_0^\infty \Big[ | w_*'|^2 + \big(\frac{4}{r^2} + \tilde f(u)\big)| w_* |^2\Big]\,r^2\,dr \ .
$$
Using decomposition $w_*(r) =u(r) \ringw (r)$ with $\ringw \in C^\infty_c (0,\infty)$, we obtain by \eqref{RS::ODE0}:
\begin{align*}
\frac{15}{32\pi} \mcQ_3(w_3)
	&= \int_0^\infty \Big\{|u\,\ringw' + u'\,\ringw|^2 + \frac{4}{r^2}\,u^2\,\ringw^2 + \tilde f(u)\,u^2\,\ringw^2\Big\}\,r^2\,dr\\
	&= \int_0^\infty \Big\{|u\,\ringw' + u'\,\ringw|^2 + \frac{4}{r^2}\,u^2\,\ringw^2 + b^2\,u^3\,\ringw^2 + f(u)\,u^2\,\ringw^2 \Big\}\,r^2\,dr\\
	&= \int_0^\infty \Big\{|u\,\ringw' + u'\,\ringw|^2 + \frac{4}{r^2}\,u^2\,\ringw^2 + b^2\,u^3\,\ringw^2 + (u'' + \frac{2}{r}u' - \frac{6}{r^2}u)\,u\,\ringw^2\Big\}\,r^2\,dr\\
	&= \int_0^\infty \Big\{|\ringw'|^2 - \frac{2}{r^2}\,\ringw^2 + b^2\,u\,\ringw^2\Big\}\,u^2\,r^2\,dr.
\end{align*}

Consider first the limit case $b^2 = 0$. \footnote{By Remark \ref{rem_scaling}, this case is equivalent to $a^2 = \infty$.} We show that $\mcQ_3(w_3)$ becomes negative by constructing a test function $w_*$ that is compactly supported close to infinity. For that, we know that  $u(r) \to s_+$ as $r \to \infty$ and therefore for small $\epsilon >0$ there exists $R> 0$ such that $s_+ (1 - \epsilon) < u(r) < s_+$ on $(R, \infty)$. The problem is solved if 
$$
\frac{15}{32\pi} \mcQ_3(w_3)\leq s_+^2 \int_R^\infty \Big\{|\ringw'|^2 - \frac{2(1-\epsilon)^2}{r^2}\,\ringw^2 \Big\} \,r^2\,dr
$$ 
is negative for  some test function $\ringw \in C^\infty_c (R, \infty)$. Since the best constant of Hardy's inequality in $\RR^3$ is $\frac{1}{4} < 2 (1-\epsilon)^2$  it is always possible to find such a function $\ringw$. \footnote{For instance, taking $n$ large enough and defining $\psi:\RR_+\to \RR_+$ by $\psi(r)=\frac{1}{R} - \frac{1}{r}$ on $(R, \frac{2nR}{n+1})$, $\psi(r)=\frac{1}{r} - \frac{1}{nR}$ on $(\frac{2nR}{n+1}, nR)$ and $\psi(r)=0$ elsewhere, we can choose $\ringw$ to be a smooth (compactly supported) approximation of $\psi$.}
Therefore $\mcQ_3 $ is not positive definite when $b^2=0$.

We turn to the case where $b^2$ is small. A perturbation argument then shows that, for any $a^2 > 0$ and $c^2 > 0$, there exists a $b_0 = b_0(a^2,c^2) > 0$ and a function $w_* \in C_c^\infty(0,\infty)$ such that $\mcQ_3(w_3)<0$ is negative for all $b^2 < b_0^2$. A rescaling argument ensures that the last statement above is equivalent to the last claim in the theorem. 
\end{proof}


\section{An orthogonal-type decomposition of $\mcQ$} 
\label{sec:SepVar}

To a certain extent, Theorem \ref{thm:scalar} suggests that in studying the stability of the second variation $\mcQ$, one cannot ignore the interaction between the various directions $E_i$. This calls for a more systematic study of the structure of $\mcQ$. As outlined in the introduction, we pursue this by performing a certain ``Fourier decomposition'' to separate variables. The reduction in the $\varphi$-variable, which is simpler, is done in subsection \ref{subsec:varphi reduction}, and the more sophisticated reduction in the $\theta$-variable is done in subsection \ref{sec-heuristics}. 

\subsection{Reduction of the $\varphi$-variable}
\label{subsec:varphi reduction}
In the previous subsection we showed that for $a^2$ large enough the radially symmetric solution $H(x)$ is unstable. Now we come back to our program of proving stability of $H$ for small $a^2$. In this section we reduce the original problem of investigating the sign of the second variation $\mcQ(V)$ to a set of simpler problems (e.g. that are independent in $\f$-variable). For this  we use a Fourier decomposition of $V$ in the $\varphi$-variable and show that the positivity of $\mcQ(V)$ can be reduced to the positivity  of some functionals $\Phi_k$ (defined in \eqref{Phidef}) for $k=0,1,2$. This type of decomposition was previously used in different contexts and for different systems, to reduce the dependence in one variable to just a few modes, see for instance \cite{Mironescu-radial}.
 
We start by the representation \eqref{defV} of $V$ as $V=\sum_{i=0}^4 w_i (r,\theta,\varphi) E_i$,
where $E_i$ are defined in \eqref{def:Ei}. Let us expand  $w_i$ using Fourier series in the $\varphi$-variable 
\begin{equation*}
w_i(r,\theta,\varphi)
	= \sum_{k=0}^\infty (\mu^{(i)}_k(r,\theta)\,\cos k\varphi + \nu^{(i)}_k(r,\theta)\,\sin k\varphi).
\end{equation*}

It is clear now that 
\be \label{Vdef}
V(r,\theta,\varphi) 
	=  \sum_{k=0}^\infty V_k(r,\theta,\varphi) = \sum_{k=0}^\infty (M_k(r,\theta, \varphi)\,\cos k\varphi + N_k(r,\theta,\varphi)\,\sin k\varphi), 
\ee
where
\be \label{Vk}
V_k(r,\theta, \varphi) = M_k(r,\theta, \varphi)\,\cos k\varphi + N_k(r,\theta,\varphi)\,\sin k\varphi,
\ee
\begin{align} \label{MN}
M_k(r,\theta,\varphi) = \sum_{i = 0}^4 \mu^{(i)}_k(r,\theta)\,E_i (\theta,\varphi)
	\,\, \text{ and } \,\,
N_k(r,\theta,\varphi) = \sum_{i = 0}^4 \nu^{(i)}_k(r,\theta)\,E_i (\theta, \varphi). 
\end{align}
We note that if $V \in H^1(\RR^3;\mcS_0)$, then $V_k \in H^1(\RR^3;\mcS_0)$ for all $k = 0, 1, \ldots$ (This can be seen by considering e.g. the contributions to $|\nabla V|^2$ in formula \eqref{import} below.)

We first prove the following lemma that reduces the study of the non-negativity of $\mcQ(V)$ to the subspace of the first three modes $V_0, V_1, V_2$.

%

\begin{lemma} \label{DVk}
Assume that, for $0 \leq k \leq 2$, $\mcQ (V_k ) \geq 0$ for any $V \in H^1(\RR^3, \mcS_0)$, where $V_k$ and $V$ are related through \eqref{Vdef}, \eqref{Vk}, \eqref{MN}. 
Then $\mcQ(V) \geq 0$ for any $V \in   H^1(\RR^3, \mcS_0)$.
\end{lemma}
\begin{proof}
Using formula \eqref{energy} for the second variation we have that 
$$\mcQ(V) = \sum_{k=0}^\infty  \mcQ(V_k).$$ 
This is because the crossing terms coming from $V_k$ and $V_j$ ($k\neq j$) involve products of $M_k$, $N_k$, $M_j$, $N_j$ that are $\f$-invariant so that
they disappear when integrating in $\f$ due to the $L^2$-orthogonality of $\{\cos k\varphi, \sin k\varphi\}$ and $\{\cos j\varphi, \sin j\varphi\}$.
After integration in the $\varphi$-variable for $k \geq 1$, we have by \eqref{energy}:
\begin{align}
\nonumber \mcQ(V_k)
	&=\pi   \int_0^\infty  \int_0^\pi \Big\{\frac{1}{3}|\partial_r \mu^{(0)}_k|^2 + |\partial_r \mu^{(1)}_k|^2 + |\partial_r \mu^{(2)}_k|^2 + |\partial_r \mu^{(3)}_k|^2 + |\partial_r \mu^{(4)}_k|^2\\
		\nonumber			&\qquad\qquad + \frac{1}{3}|\partial_r \nu^{(0)}_k|^2 + |\partial_r \nu^{(1)}_k|^2 + |\partial_r \nu^{(2)}_k|^2 + |\partial_r \nu^{(3)}_k|^2 + |\partial_r \nu^{(4)}_k|^2\\
\nonumber			&\qquad + \frac{1}{r^2}\Big[\frac{1}{3} (\partial_\theta \mu^{(0)}_k - 3 \mu^{(2)}_k)^2 + (\partial_\theta \mu^{(1)}_k - \mu^{(3)}_k)^2 + (\partial_\theta \mu^{(2)}_k + \mu^{(0)}_k - \mu^{(4)}_k)^2\\
\nonumber					&\qquad\qquad + (\partial_\theta \mu^{(3)}_k + \mu^{(1)}_k)^2 + (\partial_\theta \mu^{(4)}_k + \mu^{(2)}_k)^2\\
	\nonumber				&\qquad\qquad + \frac{1}{3} (\partial_\theta \nu^{(0)}_k - 3 \nu^{(2)}_k)^2 + (\partial_\theta \nu^{(1)}_k - \nu^{(3)}_k)^2 + (\partial_\theta \nu^{(2)}_k + \nu^{(0)}_k - \nu^{(4)}_k)^2\\
	\nonumber				&\qquad\qquad + (\partial_\theta \nu^{(3)}_k + \nu^{(1)}_k)^2 + (\partial_\theta \nu^{(4)}_k + \nu^{(2)}_k)^2 \Big]\\	
\nonumber			&\qquad + \frac{1}{r^2\,\sin^2\theta}\Big[\frac{1}{3}(k\nu^{(0)}_k + 3\sin\theta \,\mu^{(1)}_k)^2 + (k \nu^{(1)}_k - \sin\theta\,\mu^{(0)}_k - \cos\theta\,\mu^{(2)}_k - \sin \theta\,\mu^{(4)}_k)^2\\
\nonumber					&\qquad\qquad + (k \nu^{(2)}_k + \cos\theta\,\mu^{(1)}_k + \sin\theta \,\mu^{(3)}_k)^2 + (k \nu^{(3)}_k - \sin\theta\, \mu^{(2)}_k - 2\cos\theta\,\mu^{(4)}_k)^2 \\
	\nonumber				&\qquad\qquad + (k \nu^{(4)}_k + \sin\theta\,\mu^{(1)}_k + 2\cos\theta\,\mu^{(3)}_k)^2\\
\nonumber					&\qquad\qquad + \frac{1}{3}(- k \mu^{(0)}_k + 3\sin\theta \,\nu^{(1)}_k)^2 + (-k \mu^{(1)}_k - \sin\theta\,\nu^{(0)}_k - \cos\theta\,\nu^{(2)}_k - \sin \theta\,\nu^{(4)}_k)^2\\
\nonumber					&\qquad\qquad + (-k \mu^{(2)}_k + \cos\theta\,\nu^{(1)}_k + \sin\theta \,\nu^{(3)}_k)^2 + (-k \mu^{(3)}_k - \sin\theta\, \nu^{(2)}_k  - 2\cos\theta\,\nu^{(4)}_k)^2 \\
	\nonumber				&\qquad\qquad + (-k \mu^{(4)}_k + \sin\theta\,\nu^{(1)}_k + 2\cos\theta\,\nu^{(3)}_k)^2\Big]\\
	\nonumber		&\qquad + \frac{1}{3}\hat f(u)\,(|\mu^{(0)}_k|^2 + |\nu^{(0)}_k|^2) + f(u)(|\mu^{(1)}_k|^2 + |\nu^{(1)}_k|^2 + |\mu^{(2)}_k|^2 + |\nu^{(2)}_k|^2)\\
	\label{import}				&\qquad\qquad  + \tilde f(u)(|\mu^{(3)}_k|^2 + |\nu^{(3)}_k|^2 + |\mu^{(4)}_k|^2 + |\nu^{(4)}_k|^2)\Big\}\,r^2\,\sin\theta\,d\theta \,dr.
\end{align}
By the above expression, denoting
$$\tilde V_{k-1} = M_k\,\cos (k-1)\varphi + N_k\,\sin (k-1)\varphi,$$
we deduce:
\begin{align*}
&\mcQ(V_k) - \mcQ(\tilde V_{k-1})\\
	&\qquad =  \pi\int_0^\infty \int_0^{\pi} \frac{1}{\sin\theta}\Big\{(2k - 1) \Big[\frac{1}{3}(|\mu^{(0)}_k|^2 + |\nu^{(0)}_k|^2) + |\mu^{(1)}_k|^2 + |\nu^{(1)}_k|^2 + |\mu^{(2)}_k|^2 + |\nu^{(2)}_k|^2\\
					 &\qquad\qquad\qquad + |\mu^{(3)}_k|^2 + |\nu^{(3)}_k|^2 + |\mu^{(4)}_k|^2 + |\nu^{(4)}_k|^2\Big]\\
			&\qquad\qquad + 4\Big[\sin \theta(\nu^{(0)}_k \mu^{(1)}_k - \mu^{(0)}_k \nu^{(1)}_k) 
					+ \cos\theta(-\nu^{(1)}_k\,\mu^{(2)}_k + \mu^{(1)}_k\,\nu^{(2)}_k)\\
					&\qquad\qquad\qquad + \sin\theta(-\nu^{(1)}_k\,\mu^{(4)}_k + \mu^{(1)}_k\,\nu^{(4)}_k)
					+ \sin\theta( \nu^{(2)}_k\,\mu^{(3)}_k - \mu^{(2)}_k\,\nu^{(3)}_k)\\
					&\qquad\qquad\qquad + 	2\cos\theta(-\nu^{(3)}_k\,\mu^{(4)}_k + \mu^{(3)}_k\,\nu^{(4)}_k)\Big]\Big\}\,d\theta\,dr. 
\end{align*}
For $k \geq 3$, Young's inequality yields:
\begin{align*}
4\sin\theta(\nu^{(0)}_k \mu^{(1)}_k - \mu^{(0)}_k \nu^{(1)}_k) 
	&\leq \frac{2k-1}{3}  (|\mu^{(0)}_k|^2 + |\nu^{(0)}_k|^2) +  \frac{12}{2k-1}\,\sin^2\theta (|\mu^{(1)}_k|^2 + |\nu^{(1)}_k|^2)
	,\\
4 \cos\theta(-\nu^{(1)}_k\,\mu^{(2)}_k + \mu^{(1)}_k\,\nu^{(2)}_k)
	&\leq (2k - 1)\cos^2\theta (|\mu^{(1)}_k|^2 + |\nu^{(1)}_k|^2)\\
		&\qquad\qquad + (2k - 1 - \frac{4}{2k - 1})\,(|\mu^{(2)}_k|^2 + |\nu^{(2)}_k|^2)
	,\\
4 \sin\theta(-\nu^{(1)}_k\,\mu^{(4)}_k + \mu^{(1)}_k\,\nu^{(4)}_k)
	&\leq \frac{4\sin^2\theta}{2k - 1 - \frac{16}{2k-1}} (|\mu^{(1)}_k|^2 + |\nu^{(1)}_k|^2)\\
		&\qquad\qquad  + (2k - 1 - \frac{16}{2k-1})(|\mu^{(4)}_k|^2 + |\nu^{(4)}_k|^2)
		,\\
4 \sin\theta( \nu^{(2)}_k\,\mu^{(3)}_k - \mu^{(2)}_k\,\nu^{(3)}_k)
	&\leq \frac{4}{2k-1} (|\mu^{(2)}_k|^2 + |\nu^{(2)}_k|^2) + (2k - 1)\sin^2\theta (|\mu^{(3)}_k|^2 + |\nu^{(3)}_k|^2)
		,\\
8\,\cos\theta(-\nu^{(3)}_k\,\mu^{(4)}_k + \mu^{(3)}_k\,\nu^{(4)}_k)
	&\leq (2k - 1)\cos^2\theta ( |\mu^{(3)}_k|^2 + |\nu^{(3)}_k|^2 )+ \frac{16}{2k - 1}(|\mu^{(4)}_k|^2 + |\nu^{(4)}_k|^2) ,
\end{align*}
and
$$
 \frac{12}{2k-1} +  \frac{4}{2k - 1 - \frac{16}{2k-1}} < 2k-1, \, k\geq 3.
$$
Summing up the above inequalities, we get:  
$$
\mcQ(V_k) \geq  \mcQ(\tilde V_{k-1}), \, k \geq 3.
$$
By our assumption,  we know that $\mcQ(V_k) \geq 0$ for $k=0, 1, 2$ and any $V\in H^1(\RR^3, \mcS_0)$, so that we conclude that $\mcQ(V) \geq 0$.
\end{proof}

\noindent Thus, to prove stability, it remains to show that $\mcQ(V_k) \geq 0$ for $k = 0, 1, 2$ and for every $V \in H^1(\RR^3, \mcS_0)$. Below we simplify this task even further. For any $\u_0(r, \theta), \u_2(r, \theta), \u_4(r, \theta) \in H^1(\RR^3)$ that are $\f$-invariant, we define  functionals $\Phi_k$, $k=0,1,2$ (that are characteristic to the three modes $\mcQ(V_k) \geq 0$ for $k = 0, 1, 2$, as explained below):
\begin{align}
\Phi_k (\u_0, \u_2, \u_4)
	&=  \pi   \int_0^\infty  \int_0^\pi \Big\{\frac{1}{3}|\partial_r \u_0|^2   + |\partial_r \u_2|^2   +  |\partial_r \u_4|^2    + \frac{1}{r^2}\Big[\frac{1}{3} |\partial_\theta \u_0|^2 + |\partial_\theta \u_2|^2  + |\partial_\theta \u_4|^2  \Big] \nonumber \\
			&\qquad\qquad + \frac{1}{r^2}\Big[ (2 + \frac{1}{3}k^2\,\csc^2\theta) |\u_0|^2
								+ (5 + (\cot\theta + k\,\csc\theta)^2)  |\u_2|^2  \nonumber \\
						&\qquad\qquad\qquad + (2 + (2\cot\theta + k\,\csc\theta)^2)  |\u_4|^2  \Big]
						\nonumber		\\
			&\qquad\qquad + \frac{4}{r^2}\Big[   - \u_0\,( \partial_\theta  \u_2   +  \cot\theta\u_2   + k\,\csc \theta\, \u_2  ) \nonumber \\
							&\qquad\qquad\qquad   +  ( - \partial_\theta  \u_2  + \cot\theta\, \u_2 + k\,\csc\theta\, \u_2  ) \, \u_4   	\Big] \nonumber \\
  			&\qquad\qquad + \frac{1}{3}\hat f(u)\, |\u_0|^2  + f(u)  |\u_2|^2 
					  + \tilde f(u)  |\u_4|^2  \Big\}\,r^2\,\sin\theta\,d\theta \,dr.
					   \label{Phidef}
\end{align}
It should be noted that, in \eqref{Phidef}, we use the convention that
\[
\Phi_k (\u_0, \u_2, \u_4) = +\infty \,\,\text{ if }\,\, \frac{1}{r}k\csc\theta|\u_0| + \frac{1}{r}|\cot\theta + k\,\csc\theta|\,|\u_2|  + \frac{1}{r}|2\cot\theta + k\,\csc\theta|\,|\u_4| \notin L^2(\RR^3).
\]
(This is justified by the observation that the expression inside the curly bracket in the integral on the right hand side of \eqref{Phidef} is bounded from below by
\[
\frac{k^2\csc^2\theta}{3r^2}\,|\u_0|^2 + \frac{(\cot\theta + k\,\csc\theta)^2}{2r^2}|\u_2|^2 + \frac{(2\cot\theta + k\,\csc\theta)^2}{2r^2}\,|\u_4|^2 - C(1 + r^{-2}) [|\u_0|^2 + |\u_2|^2 + |\u_4|^2].)
\]

To conclude this subsection, we establish:
\begin{proposition} \label{Phi012}
The second variation $\mcQ $ is nonnegative definite if and only if the functionals $\Phi_0$, $\Phi_1$ and $\Phi_2$ are nonnegative. 
\end{proposition} 
\begin{proof} By \eqref{import}, we observe that $\mcQ(V_k)$ splits in the following way:
\begin{align*}
\mcQ(V_k) 
	&= \mcQ((\mu^{(0)}_k E_0+ \mu^{(2)}_k E_2 + \mu^{(4)}_k E_4) \cos k \varphi + (\nu^{(1)}_k E_1 + \nu^{(3)}_k E_3) \sin k \varphi ) \\
	&+ \mcQ ((-\nu^{(0)}_k E_0- \nu^{(2)}_k E_2 - \nu^{(4)}_k E_4) \cos k \varphi + (\mu^{(1)}_k E_1 + \mu^{(3)}_k E_3) \sin k \varphi),
\end{align*}
where we recall that $\mu_k^{(j)}$ and $\nu_k^{(j)}$ are $\f$-independent. It is thus sufficient to investigate the sign of the functionals
\begin{align}
\nonumber \mcQ_k
	&:=\mcQ_k(\xi_0, \xi_1, \xi_2, \xi_3, \xi_4) = \mcQ((\xi_0 E_0 +  \xi_2 E_2 + \xi_4 E_4)\cos k\varphi + (\xi_1 E_1 + \xi_3 E_3)\sin k \varphi)\\
	&=\pi   \int_0^\infty  \int_0^\pi \Big\{\frac{1}{3}|\partial_r \xi_0|^2 + |\partial_r \xi_1|^2 + |\partial_r \xi_2|^2 + |\partial_r \xi_3|^2 + |\partial_r \xi_4|^2\nonumber\\
	\nonumber			&\qquad + \frac{1}{r^2}\Big[\frac{1}{3} (\partial_\theta \xi_0 - 3 \xi_2)^2 + (\partial_\theta \xi_2 + \xi_0 - \xi_4)^2\\
\nonumber					&\qquad\qquad + (\partial_\theta\xi_4 + \xi_2)^2 + (\partial_\theta \xi_1 - \xi_3)^2+ (\partial_\theta \xi_3 + \xi_1)^2 \Big]\\
\nonumber			&\qquad + \frac{1}{r^2\,\sin^2\theta}\Big[(k\xi_1 - \sin\theta\,\xi_0 - \cos\theta\,\xi_2 - \sin \theta\,\xi_4)^2\\
\nonumber					&\qquad\qquad +  (k \xi_3 - \sin\theta\, \xi_2 - 2\cos\theta\,\xi_4)^2 \\
	\nonumber					&\qquad\qquad + \frac{1}{3}(- k \xi_0 + 3\sin\theta \,\xi_1)^2 +(-k \xi_2 + \cos\theta\,\xi_1 + \sin\theta \,\xi_3)^2  \\
	\nonumber				&\qquad\qquad + (-k \xi_4 + \sin\theta\,\xi_1 + 2\cos\theta\,\xi_3)^2\Big]\\
	\nonumber		&\qquad + \frac{1}{3}\hat f(u)\,\xi_0^2 + f(u)(\xi_1^2 + \xi_2^2) + \tilde f(u)(\xi_3^2 + \xi_4^2) \Big\}\,r^2\,\sin\theta\,d\theta \,dr,
\end{align}
where $\{ \xi_i=\xi_i(r, \theta) \}_{i=0}^4\subset H^1(\RR^3)$ and $\mcQ_k$ can take value $+\infty$.

We define the following change of variables \footnote{As we will see later, this change of variables is done in order that a certain system of ODEs
(related to the Euler-Lagrange equations associated to the functional $\Phi_k$) decouples. See Remark \ref{rem:17VI13.01}.}
\begin{equation}
\label{chg_var}
\u_0 = \frac{1}{2}\xi_0, \u_1 = \frac{1}{2}(\xi_1 + \xi_2), \u_2 = \frac{1}{2}(\xi_1 - \xi_2), \u_3 = \frac{1}{2}(\xi_3 + \xi_4) \text{ and } \u_4 = \frac{1}{2}(\xi_3 - \xi_4).
\end{equation}
By the above change of variable, we obtain:
\begin{align*}
\mcQ_k 
	&  =\pi   \int_0^\infty  \int_0^\pi \Big\{\frac{4}{3}|\partial_r \u_0|^2   + 2|\partial_r \u_1|^2    + 2|\partial_r \u_2|^2 +  2|\partial_r \u_3|^2   + 2|\partial_r \u_4|^2  \\
			&\qquad\qquad + \frac{1}{r^2}\Big[\frac{4}{3} |\partial_\theta \u_0|^2 + 2|\partial_\theta \u_1|^2 + 2 |\partial_\theta \u_2|^2 + 2|\partial_\theta \u_3|^2 + 2 |\partial_\theta \u_4|^2\Big]\\
			&\qquad\qquad + \frac{1}{r^2}\Big[ 4(2 + \frac{1}{3}k^2\,\csc^2\theta) |\u_0|^2
								+ 2(4 + (k^2 + 1)\,\csc^2\theta) (|\u_1|^2 + |\u_2|^2)\\
						&\qquad\qquad\qquad + 2(-2 + (k^2 + 4)\,\csc^2\theta) (|\u_3|^2 + |\u_4|^2)\Big]
								\\
			&\qquad\qquad + \frac{4}{r^2}\Big[ 2\u_0\,(\partial_\theta (\u_1 - \u_2) + \cot\theta\,(\u_1 - \u_2) - k\,\csc \theta\,(\u_1 + \u_2))\\
							&\qquad\qquad\qquad   +  2( - \partial_\theta  \u_1  + \cot\theta\, \u_1 - k\,\csc\theta\, \u_1  ) \, \u_3 \\
							&\qquad\qquad\qquad  +  2(- \partial_\theta  \u_2  + \cot\theta\, \u_2  + k\,\csc\theta\, \u_2 ) \, \u_4   	\Big]
\end{align*}
\begin{align*}
			&\qquad\qquad\qquad + \frac{2}{r^2}\Big[  
						- 2k\cot\theta\,\csc\theta\, (|\u_1|^2 - |\u_2|^2)
						   - 4k \cot\theta\,\csc\theta \,(|\u_3|^2 - |\u_4|^2)\Big]\\
 			&\qquad\qquad\qquad + \frac{4}{3}\hat f(u)\, |\u_0|^2  + 2f(u)( |\u_1|^2 + |\u_2|^2 ) 
					  + 2\tilde f(u)(  |\u_3|^2 + |\u_4|^2  )\Big\}\,r^2\,\sin\theta\,d\theta\, dr,
\end{align*}
where we used integration by parts for the terms $\partial_\theta \u_0 (\u_2-\u_1) \sin \theta$, $\partial_\theta \u_3 \u_1 \sin \theta$ and
$\partial_\theta \u_4 \u_2 \sin \theta$. 
After rearranging the terms we obtain
$$
\mcQ_k = 2\Phi_k (\tilde \u_0,\tilde \u_1,\tilde \u_3) + 2\Phi_k(\u_0, \u_2, \u_4)
$$
where $\Phi_k(\u_0, \u_2, \u_4)$ is defined in \eqref{Phidef} and $\tilde \u_0 (r,\theta)= \u_0 (r, \pi - \theta)$, $\tilde \u_1(r,\theta)= \u_1 (r, \pi - \theta)$, $\tilde \u_3 (r,\theta)= - \u_3 (r, \pi - \theta)$.
\end{proof}
\subsection{Separating $r$ and $\theta$ variables}
\label{sec-heuristics}

We continue with the analysis of the sign of $\Phi_k(\u_0, \u_2, \u_4)$ ($k=0,1,2$) where the functions $\u_0, \u_2, \u_4$ are $\f$-independent. The idea is to separate variables in $\u_m (r, \theta)$, $m=0, 2, 4$. Before giving a precise statement, let us explain the ideas behind the construction. 

Heuristically speaking, we would like to ``minimize'' $\Phi_k$. Since $\Phi_k$ is quadratic, we need to subject it to some normalization. It is instructive to notice that the coefficients in front of $\u_m^2$ in the expression for $\Phi_k$ behave differently at infinity for different $m$. This suggests that a natural normalization is
\[
\pi\int_0^\infty\int_0^\pi \Big[\u_0^2 + \frac{1}{r^2}\u_2^2 + \u_4^2\Big] r^2\,\sin\theta\, d\theta\,dr = 1.
\]
The \emph{formal} Euler-Lagrange equations associated to such minimization problem read:

\begin{align*}
&-\frac{1}{3r^2}\partial_r(r^2\,\partial_r \u_0)
	- \frac{1}{3r^2\sin\theta}\partial_\theta(\sin\theta \partial_\theta \u_0)\\
	&\qquad + \frac{1}{r^2}(2+\frac 1 3 k^2 \csc^2 \theta) \u_0 
		- \frac{2}{r^2}(\partial_\theta \u_2 +\cot \theta \u_2 + k \csc \theta \u_2)\\
	&\qquad + \frac{1}{3}\hat f(u)\,\u_0
		= \lambda \u_0,
\end{align*}

\begin{align*}
&-\frac{1}{r^2}\partial_r (r^2\,\partial_r \u_2)
	- \frac{1}{r^2\sin\theta}\partial_\theta(\sin\theta \partial_\theta \u_2)\\
	&\qquad + \frac{1}{r^2}[5+(\cot \theta+k\csc \theta)^2] \u_2
		\\
	&\qquad	+ \frac{2}{r^2\sin\theta}[\partial_\theta (\sin \theta \u_0)+\partial_\theta (\sin \theta \u_4) +(\cos\theta+k)(- \u_0 + \u_4)]\\
	&\qquad + f(u)\,\u_2
		= \frac{\lambda}{r^2}\,\u_2
	,
\end{align*}

\begin{align*}
&-\frac{1}{r^2}\partial_r(r^2\,\partial_r \u_4)
	- \frac{1}{r^2\sin\theta}\partial_\theta(\sin\theta \partial_\theta \u_4)\\
	&\qquad + \frac{1}{r^2}[2+(2\cot \theta+k \csc\theta)^2] \u_4 
		+\frac{2}{r^2}(-\partial_\theta \u_2 +\cot \theta \u_2 + k \csc \theta \u_2)\\
	&\qquad + \tilde f(u)\,\u_4
		= \lambda \u_4.
\end{align*}

To proceed with our heuristic discussion, we formally solve these equations at infinity by making an ansatz that $\u_m (r,\theta) \sim \frac{1}{r^{\alpha_m}}\, u_m(\theta)$ as $r \to \infty$ with $m=0,2,4$. 
Keeping only leading terms as $r \rightarrow \infty$ in the above system, we obtain the relations 
$$
\alpha_0  = \alpha_4 = \alpha_2 + 2,
$$
and
\begin{align}
&u_0 \sim  \partial_\theta u_2 + \cot\theta u_2 + k\,\csc\theta\,u_2
	,\label{Eq:u0ODE}\\
&u_4 \sim -\partial_\theta u_2 + \cot\theta u_2 + k\,\csc\theta\,u_2
	,\label{Eq:u4ODE}\\
&u_2 \sim -\frac{1}{\sin\theta}\partial_\theta(\sin\theta \partial_\theta u_2) + (\cot\theta + k\csc\theta)^2 u_2 .
\label{Eq:u2ODE}
\end{align}

Therefore, for $k=0,1,2$, we consider the spectral problem for the operator
\begin{equation}
\label{EVP:U2}
T^{(2)}_{k} \equiv  -\frac{1}{\sin\theta}\partial_\theta(\sin\theta \partial_\theta  ) + [1+(\cot\theta + k\csc\theta)^2],
\end{equation}
or more precisely, on the couples (eigenfunction, eigenvalue)$=(u_{k,i}^{(2)}, \lambda_{k,i})$
$$
T^{(2)}_{k} u_{k,i}^{(2)} = \lambda_{k,i}  u_{k,i}^{(2)}, \quad i\geq 1.
$$
For each $k=0,1,2$, these eigenfunctions $\{ u_{k,i}^{(2)} \}_{i\geq 1}$ will form a basis in  $L^2((0,\pi);\sin\theta\,d\theta)$ and therefore we have the following expansion for $\u_2$
$$
\u_2 (r,\theta) = \sum_{i\geq 1} w_{k,i}^{(2)} (r) u_{k,i}^{(2)} (\theta).
$$
The key observation is the following: For $k=0,1,2$, if $u^{(2)}_{k,i}$ is an eigenfunction of the operator $T^{(2)}_{k}$ associated to an eigenvalue $\lambda_{k,i}$, and if (compare \eqref{Eq:u0ODE}, \eqref{Eq:u4ODE})
\begin{align}
&u^{(0)}_{k,i} :=  \partial_\theta u^{(2)}_{k,i} + \cot\theta u^{(2)}_{k,i} + k\,\csc\theta\,u^{(2)}_{k,i}
	,\label{Def:U0}\\
&u^{(4)}_{k,i} = -\partial_\theta u^{(2)}_{k,i} + \cot\theta u^{(2)}_{k,i} + k\,\csc\theta\,u^{(2)}_{k,i}
	,\label{Def:U4}
\end{align}
then $u^{(0)}_{k,i}$ and $u^{(4)}_{k,i}$ are eigenfunctions of the following operators
\begin{align}
T^{(0)}_{k} \equiv &-\frac{1}{\sin\theta}\partial_\theta(\sin\theta \partial_\theta ) +  k^2\csc^2\theta \, 
	,\label{EVP:U0}\\
T^{(4)}_{k} \equiv &-\frac{1}{\sin\theta}\partial_\theta(\sin\theta \partial_\theta ) +  [4 + (2\cot\theta + k\csc\theta)^2 ]\, ,
\label{EVP:U4}
\end{align}
and satisfy
\be
\label{eigenvalu}
T^{(0)}_{k} u_{k,i}^{(0)} =\lambda_{k,i} u_{k,i}^{(0)}, \quad T^{(4)}_{k} u_{k,i}^{(4)} =\lambda_{k,i} u_{k,i}^{(4)}.
\ee
This is a lengthy but easy computation, which can be observed heuristically by noting their compatibility with the the second order leading terms as $r\to \infty$ in the above mentioned formal Euler-Lagrange equations associated to $\Phi_k$.

The functions $\{u^{(0)}_{k,i}\}$ and $\{u^{(4)}_{k,i}\}$ also form bases of $L^2((0,\pi);\sin\theta\,d\theta)$ and therefore one can write $\u_m$ as
\[
\u_m(r,\theta) = \sum_i w^{(m)}_{k,i}(r)\,u^{(m)}_{k,i}(\theta), \quad m=0,4.
\]
A remarkable consequence of this decomposition is that the functionals $\Phi_k$'s also decouple as
\[
\Phi_k (\u_0, \u_2, \u_4) = \pi \sum_{i \geq C(k)} \Phi_{0,i} (w^{(0)}_{k,i}, w^{(2)}_{k,i}, w^{(4)}_{k,i}), \quad k=0,1,2,
\]
where functionals $\Phi_{0,i}$ are defined in \eqref{new_form_phi} below and $C(k)$ is a natural number depending on $k$. 

This decoupling is unexpected, because we deal with a system and it is a priori not clear how the different bases for different $\u_m$'s interact with one another.

\begin{remark}\label{rem:17VI13.01}
We would like to note that the change of variables \eqref{chg_var} was done in order that the above system decouples. If one attempts to directly apply the above procedure to $\mcQ_k$, it does not seem that one would obtain the above separation of variables. For example, it is not clear that the equations of the relevant $\xi_k$'s would decouple as in \eqref{Eq:u0ODE}-\eqref{Eq:u2ODE}.
\end{remark}

We now proceed to formalize the above heuristic argument. As mentioned above, the functions $u^{(2)}_{k,i}$ are eigenfunctions of the operators $T_k^{(2)}$ defined in \eqref{EVP:U2}. To obtain a well-posed eigenproblem, we need to furnish a boundary condition. Note that the operators $T_k^{(2)}$ are Fuchsian at both $\theta = 0$ and $\theta = \pi$. We thus impose boundary conditions so that the less regular Fuchsian indices drop out. More precisely, we define the following spectral problems for $u^{(2)}_{k,i}$, $i \geq 1$:
\begin{equation}
\left\{\begin{array}{l}
\displaystyle -\frac{1}{\sin\theta}\partial_\theta(\sin\theta \partial_\theta u^{(2)}_{0,i}) + \csc^2\theta\,u^{(2)}_{0,i} = \lambda_{0,i}\,u^{(2)}_{0,i}
	,\\
\displaystyle u^{(2)}_{0,i}(0) = u^{(2)}_{0,i}(\pi) = 0, \qquad \int_0^\pi |u^{(2)}_{0,i}|^2\,\sin\theta\,d\theta = 1.
\end{array}\right.
	\label{EVP:k0U2}
\end{equation}

\begin{equation}
\left\{\begin{array}{l}
\displaystyle -\frac{1}{\sin\theta}\partial_\theta(\sin\theta \partial_\theta u^{(2)}_{1,i}) + (1 + (\cot\theta + \csc\theta)^2)\,u^{(2)}_{1,i} = \lambda_{1,i}\,u^{(2)}_{1,i}
	,\\
\displaystyle \partial_\theta u^{(2)}_{1,i}(0) = \partial_\theta  u^{(2)}_{1,i}(\pi) = 0, \qquad \int_0^\pi |u^{(2)}_{1,i}|^2\,\sin\theta\,d\theta = 1.
\end{array}\right.
	\label{EVP:k1U2}
\end{equation} 
\begin{equation}
\left\{\begin{array}{l}
\displaystyle -\frac{1}{\sin\theta}\partial_\theta(\sin\theta \partial_\theta u^{(2)}_{2,i}) + (1 + (\cot\theta + 2\csc\theta)^2)\,u^{(2)}_{2,i} = \lambda_{2,i}\,u^{(2)}_{2,i}
	,\\
\displaystyle   u^{(2)}_{2,i}(0) =   u^{(2)}_{2,i}(\pi) = 0, \qquad \int_0^\pi |u^{(2)}_{2,i}|^2\,\sin\theta\,d\theta = 1.
\end{array}\right.
	\label{EVP:k2U2}
\end{equation}
The functions $u^{(0)}_{k,i}$ and $u^{(4)}_{k,i}$ are defined as in \eqref{Def:U0} and \eqref{Def:U4}.

\begin{proposition}\label{prop:Phikdecomp}
Let $\lambda_{k,i}$, $u^{(0)}_{k,i}$, $u^{(2)}_{k,i}$ and $u^{(4)}_{k,i}$ be defined as above. The following assertions are true.

(a) $\lambda_{0,i} = i(i+1)$,  and the set $\{ u^{(2)}_{0,i} \}_{i=1}^\infty$ is an orthonormal basis in $L^2((0,\pi);\sin\theta\,d\theta)$. Moreover, the sets $\{ u^{(0)}_{0,i} \}_{i=1}^\infty \cup \{u^{(0)}_{0,0}\equiv 1\}$ and $\{ u^{(4)}_{0,i} \}_{i=2}^\infty$ are orthogonal bases of $L^2((0,\pi);\sin\theta\,d\theta)$.

\smallskip
(b) $\lambda_{1,i} = i(i+1)$ and $\{ u^{(2)}_{1,i} \}_{i=1}^\infty$, $\{ u^{(0)}_{1,i} \}_{i=1}^\infty$ and $\{ u^{(4)}_{1,i} \}_{i=2}^\infty$ are orthogonal bases of $L^2((0,\pi);\sin\theta\,d\theta)$.

\smallskip
(c) $\lambda_{2,i} =(i+1)(i+2)$, and $\{ u^{(2)}_{2,i} \}_{i=1}^\infty$, $\{ u^{(0)}_{2,i} \}_{i=1}^\infty$ and $\{ u^{(4)}_{2,i} \}_{i=1}^\infty$ are orthogonal bases of $L^2((0,\pi);\sin\theta\,d\theta)$.

\smallskip
(d)  Fix $k \in \{0,1,2\}$. Let $\u_0, \u_2, \u_4 \in H^1(\RR^3)$ be $\varphi$-independent functions satisfying 
\[
\frac{1}{r}k\csc\theta|\u_0| + \frac{1}{r}|\cot\theta + k\,\csc\theta|\,|\u_2|  + \frac{1}{r}|2\cot\theta + k\,\csc\theta|\,|\u_4| \in L^2(\RR^3).
\]
Then
$$
\u_m (r, \theta) = \sum_i w_{k,i}^{(m)} (r) u_{k,i}^{(m)} (\theta),
$$
where $w_{k,i}^{(m)} \in H^1((0,\infty);r^2\,dr)$ and
\begin{align*}
&\sum_i c_{k,i}^{(m)}\Big[ \|w_{k,i}^{(m)}\|_{H^1((0,\infty);r^2\,dr)}^2 + \lambda_{k,i} \|\frac{1}{r}w_{k,i}^{(m)}\|_{L^2((0,\infty);r^2dr)}^2\Big]  < \infty,
\end{align*}
where
\[
c_{k,i}^{(m)} = \left\{\begin{array}{ll}
\lambda_{k,i}& \text{ if } m = 0,\\
1 & \text{ if } m = 2,\\
\lambda_{k,i} - 2 & \text{ if } m = 4.
\end{array}\right.
\]

\noindent The functionals $\Phi_k$ can then be written as
\begin{align*}
 \Phi_0(\u_0, \u_2, \u_4) 
 	&= \pi \sum_{i=0}^\infty \Phi_{0,i}(w^{(0)}_{0,i}, w^{(2)}_{0,i}, w^{(4)}_{0,i}),\\
  \Phi_1(\u_0, \u_2, \u_4) 
  	&= \pi \sum_{i=1}^\infty \Phi_{0,i}(w^{(0)}_{1,i}, w^{(2)}_{1,i}, w^{(4)}_{1,i}),\\
\Phi_2(\u_0, \u_2, \u_4) 
	&=  \pi \sum_{i=1}^\infty \Phi_{0,i+1}(w^{(0)}_{2,i}, w^{(2)}_{2,i}, w^{(4)}_{2,i}),
\end{align*}
where
\begin{align}
\Phi_{0,0} (w_0, w_2, w_4) 
	&= \frac{2}{3}\int_0^\infty \Big[|\partial_r w_0|^2 + \frac{6}{r^2}\,|w_0|^2 + \hat f(u)\,|w_0|^2\Big]\,r^2\,dr
	, \nonumber \\
\Phi_{0,i} (w_0, w_2, w_4) 
	&=  \int_0^\infty    \Big\{\frac{\lambda_{0,i}}{3}|\partial_r w_0|^2   + |\partial_r w_2|^2   +  (\lambda_{0,i} - 2)|\partial_r w_4|^2    \nonumber \\
			&\qquad + \frac{1}{r^2}\Big[\frac{\lambda_{0,i}(\lambda_{0,i} + 6)}{3} |w_0|^2 + (\lambda_{0,i} + 4) |w_2|^2  + (\lambda_{0,i} - 2)^2 \,|w_4|^2  
							\nonumber	\\
					&\qquad\qquad\qquad   - 4\lambda_{0,i}\,w_0\,w_2 + 4(\lambda_{0,i} - 2) \,w_2\,w_4 \Big] \nonumber \\
  			&\qquad + \frac{\lambda_{0,i}}{3}\hat f(u)\, |w_0|^2  + f(u)  |w_2|^2 + (\lambda_{0,i}  - 2)\tilde f(u)  |w_4|^2  \Big\}\,r^2 \,dr, \label{new_form_phi}
\end{align}
for $i\geq 1$ where $w_0, w_2, w_4$ depend only on $r$-variable and belong to $H^1(\RR^3)$.
\end{proposition}

The proof of Proposition \ref{prop:Phikdecomp} is technical and we postpone it to Section \ref{sec:bases}.

\begin{remark}
We will see in Proposition \ref{prop:phi0,2reduction} that $\Phi_{0,i} \geq 0$ for $i \notin \{2,3\}$. In particular, the infinite sums $\sum \Phi_{0,i}$ in Proposition \ref{prop:Phikdecomp} make sense.
\end{remark}

\section{On the non-negativity of the reduced functionals}\label{sec:stability}
In the previous section we showed that in order to prove non-negativity of the second variation $\mcQ(V)$ it is enough to prove non-negativity of the functionals $\Phi_{0,i}$ defined in \eqref{new_form_phi} for $i\geq 0$.

\subsection{Positivity of $\Phi_{0,i}$ for $i \notin \{2, 3\}$}
\label{subsec:phi0,i}

In this section we exclude the cases $i=2,3$ and prove the following proposition:
\begin{proposition}
\label{prop:phi0,2reduction}
Assume $i \in \NN$ and $i\not\in\{2,3\}$. Then
\[
\Phi_{0,i} (w_0, w_2, w_4) \geq 0 \text{ for all } w_0, w_2, w_4 \in H^1((0,\infty);r^2\,dr).
\]
Moreover $\Phi_{0,i} (w_0, w_2, w_4)=0$ if and only if $i=1$ and $w_0(r)= Au'$, $w_2(r)=\frac{2Au}{r}$, where $A \in \RR$ and $u$ solves \eqref{RS::ODE0}, \eqref{BC}.  
\end{proposition}

\begin{proof} Note that $C_c^\infty(0,\infty)$ is dense in $H^1((0,\infty);r^2\,dr)$. Thus in proving the non-negativity of $\Phi_{0,i}$, it suffices to consider $w_0, w_2, w_4 \in C_c^\infty(0,\infty)$.

In subsection~\ref{instab} we  already showed that $\Phi_{0,0}(w_0) > 0$ for all $w_0 \in C_c^\infty(0,\infty)$, $w_0 \not\equiv 0$ (see the proof of \eqref{mcE0positive}).

Let us show that $\Phi_{0,1} $ is non-negative. Using proposition~\ref{prop:Phikdecomp} we have $\lambda_{0,1}=2$. 
Taking any $w_0, w_2 \in C_c^\infty(0,\infty)$ and using Hardy decomposition (see Section~\ref{instab} or Appendix, Lemma~\ref{thm:hardy}) $w_0=u' \xi$ and $w_2=u \eta$, where $\xi, \eta \in C_c^\infty(0,\infty)$ and $u$ solves \eqref{RS::ODE0}, \eqref{BC}, we obtain
\begin{align*}
\Phi_{0,1} (w_0, w_2) =\frac{1}{3} \int_0^\infty \Big\{2|u'|^2\,|\xi'|^2 + 3|u|^2\,|\eta'|^2 + \frac{24}{r^3}u\,u'\,|\xi|^2 - \frac{24}{r^2}\,u\,u'\,\xi\,\eta \Big\}\,r^2\,dr.
\end{align*}
It is clear that $-\frac{24}{r^2}\,u\,u'\,\xi\,\eta \geq - \frac{24}{r}u\,u'\,|\frac{\xi}{r}|^2 - \frac{6}{r}u\,u'\,|\eta|^2$ with equality achieved only for $\eta = {2\xi \over r}$. Therefore
$$
\Phi_{0,1} (w_0, w_2) \geq \frac{1}{3}\int_0^\infty \Big\{2|u'|^2\,|\xi'|^2 + 3|u|^2\,|\eta'|^2  - \frac{6}{r}u\,u'\,|\eta|^2   \Big\}\,r^2\,dr.
$$
Integration by parts yields
$$
- \int_0^\infty 6{r}u\,u'\,|\eta|^2\, dr = \int_0^\infty \Big\{   3 u^2 \eta^2 + 6 u^2\,r \eta \eta'    \Big\} \, dr,
$$
and therefore
$$
\Phi_{0,1} (w_0, w_2) \geq \frac{1}{3}\int_0^\infty \Big\{2|u'|^2\,|\xi'|^2 + 3u^2 \left( \eta + r \eta'  \right)^2   \Big\}\,r^2\,dr \geq 0.
$$
Note that the above estimates hold for $w_0, w_2 \in H^1((0,\infty);r^2\,dr)$ by a standard density argument. 
It can then be checked that $\Phi_{0,1} (w_0, w_2) = 0$ if and only if $w_0(r)= Au'$, $w_2(r)=\frac{2Au}{r}$ for some $A \in \RR$.

We will now show that $\Phi_{0,i}$ is positive definite for $i \geq 4$ (or $\lambda_{0,i} \geq 20$). We note the following inequality, which holds for $\lambda \geq 20$:
\begin{align*}
&\frac{1}{3}\lambda^2 x^2 + (\lambda - 2)\,y^2 + (\lambda - 2)(\lambda - 8)\,z^2
	 + 4\lambda\,xy + 4(\lambda - 2)\,yz\\
	&\qquad \geq \frac{1}{3}\lambda^2 x^2 + (\lambda - 2)[1 - \frac{4}{\lambda - 8}]\,y^2
	 + 4\lambda\,xy
	\geq 0.
\end{align*}
It follows that, when $\lambda_{0,i} \geq 20$,
\begin{align*}
\Phi_{0,i}(w_0, w_2, w_4)
	&\geq \int_0^\infty \Big\{\frac{\lambda_{0,i}}{3}|\partial_r w_0|^2 + |\partial_r w_2|^2 + (\lambda_{0,i} - 2)|\partial_r w_4|^2\\
		&\qquad+ \frac{1}{r^2}\Big[2\lambda_{0,i}\,|w_0|^2 + 6\,|w_2|^2 + 6(\lambda_{0,i} - 2)\,|w_4|^2\Big]\\
		&\qquad + \frac{\lambda_{0,i}}{3}\,\hat f(u)\,|w_0|^2 + f(u)\,|w_2|^2 + (\lambda_{0,i} - 2)\tilde f(u)\,|w_4|^2\Big\}\,r^2\,dr\\
	&=\left[ \frac{\lambda_{0,i}}{2}\Phi_{0,0}(w_0) + \Phi_{0,1}(0, w_2) + (\lambda_{0,i} - 2) \Phi_{0,1}(0, w_4)\right] \\ 
		&\qquad + \int_0^\infty (\lambda_{0,i} - 2)\,b^2\,u\,|w_4|^2\,r^2\,dr \geq 0
	.
\end{align*}
Equality holds if and only if  $w_0 = w_2 = w_4 \equiv 0$. 
\end{proof}


\noindent We have shown that $\Phi_{0,i}$ is non-negative definite for $i \notin \{ 2, 3 \}$. It remains to consider functionals $\Phi_{0,3}$ and $\Phi_{0,2}$  that we do in the next section.


\subsection{Positivity of $\Phi_{0,3}$ and $\Phi_{0,2}$}
\label{subsec:phi0,2}

In this subsection we consider the functionals $\Phi_{0,3}$ and $\Phi_{0,2}$. We first show that the positivity of $\Phi_{0,2}$ implies that of $\Phi_{0,3}$ and subsequently show that $\Phi_{0,2}$ is indeed positive.
Recall that
\begin{align*}
\Phi_{0,3}(w_0,w_2,w_4) 
	&= \int_0^\infty \Big\{4|\partial_r w_0|^2 + |\partial_r w_2|^2 + 10|\partial_r w_4|^2\\
		&\qquad + \frac{1}{r^2}\Big[72|w_0|^2 + 16|w_2|^2 + 100|w_4|^2\\	
			&\qquad\qquad\qquad - 48 w_0 w_2 + 40 w_2 w_4\Big]\\
		&\qquad + 4\hat f(u) |w_0|^2 + f(u)| w_2|^2 + 10\tilde f(u)|w_4|^2\Big\}r^2\,dr,\\
\Phi_{0,2}(w_0, w_2, w_4)
	&= \int_0^\infty \Big\{2|\partial_r w_0|^2 + |\partial_r w_2|^2 + 4|\partial_r w_4|^2\\
		&\qquad + \frac{1}{r^2}\Big[24|w_0|^2 + 10|w_2|^2 + 16|w_4|^2\\	
			&\qquad\qquad\qquad - 24 w_0 w_2 + 16 w_2 w_4\Big]\\
		&\qquad + 2\hat f(u) |w_0|^2 + f(u)| w_2|^2 + 4\tilde f(u)|w_4|^2\Big\}r^2\,dr.
\end{align*}

\begin{lemma} \label{Phi3,2}
For any $w_0, w_2, w_4 \in H^1( (0, \infty); r^2\,dr) $ the following inequality holds
\[
\Phi_{0,3} (w_0, w_2, w_4) \geq \Phi_{0,2} (\sqrt{2}w_0, w_2, \frac{\sqrt{10}}{2}w_4) .
\]
\end{lemma}
\begin{proof}  Defining $\tilde w_0 = \sqrt{2}w_0$ and $\tilde w_4 = \frac{\sqrt{10}}{2}w_4$ we have:
\begin{align*}
&\Phi_{0,3}(w_0, w_2, w_4) - \Phi_{0,2}(\tilde w_0, w_2, \tilde w_4)\\
	&\qquad =  \int_0^\infty \Big[12|\tilde w_0|^2 + 6|w_2|^2 + 24|\tilde w_4|^2
			- 24(\sqrt{2} - 1) \tilde w_0 w_2 + 8(\sqrt{10} - 2) w_2 \tilde w_4\Big]\,dr\\
	&\qquad  \geq \int_0^\infty \Big(6 - 12(\sqrt{2} - 1)^2 - \frac{2(\sqrt{10} - 2)^2}{3}\Big)\int_0^\infty |w_2|^2\,dr
		\geq 0.
\end{align*}
\end{proof}

The main result of this subsection is formulated in the following proposition:
\begin{proposition}\label{prop:phi2coercivity} For any $b^2,c^2>0$ there exist $\bar a^2 >0$ and $\delta_0 > 0$ such that for $a^2<\bar a^2$ and for all $w_0, w_2, w_4 \in H^1((0,\infty);r^2\,dr)$ the following inequality holds
$$
\Phi_{0,2} (w_0, w_2, w_4) \geq \delta_0 \int_0^\infty (w_0^2 + w_2^2+w_4^2)\, dr.
$$ 
\end{proposition}
The proof of this proposition relies on the following three lemmas that we prove below.
\begin{lemma}\label{lemma:phi22}
For any $w_2 \in C_c^\infty(0,\infty)$ the following inequality holds
\[
\Phi_{0,2} (0, w_2, 0) \geq \int_0^\infty  \left(4\frac{u' r}{u} + 2 \right)w_2^2\, dr
\]
\end{lemma}
\begin{proof}
We want to use a special Hardy decomposition to prove lemma. Let us define $v=\frac{u}{r^\alpha}$ and see what equation is satisfied by $v$. Here $u(r)$ is the solution of \eqref{RS::ODE0},  and $\alpha \in \RR$ to be specified later. It is straightforward to compute
$$
v'=\frac{u'}{r^\alpha} - \frac{\alpha u}{r^{\alpha+1}}, \quad v'' = \frac{u''}{r^\alpha} - \frac{2 \alpha u'}{r^{\alpha+1}}+  \frac{\alpha(\alpha+1) u}{r^{\alpha+2}}.
$$
From equation \eqref{RS::ODE0} for $u$ we know that
$$
u''=-\frac{2 u'}{r}+\frac{6u}{r^2}+f(u)u.
$$
Therefore after a simple calculation we obtain
$$
v'' + \frac{2(1+ \alpha)v'}{r} - \frac{(6-\alpha(\alpha+1))v}{r^2} = f(u)v
$$
Take any $w_2 \in  C_c^\infty(0,\infty)$ and use the  Hardy decomposition $w_2 =v \eta$, where $\eta \in C_c^\infty(0, \infty)$. It is straightforward to obtain
\begin{align*}
\Phi_{0,2}(0,w_2,0) &=\int_0^\infty \left( |w_2'|^2\, r^2 + 10w_2^2+f(u) w_2^2 \, r^2 \right) \, dr  \nonumber \\
	&=\int_0^\infty \Big[ |v' \eta + v \eta'|^2 r^2 + 10v^2 \eta^2 \nonumber\\ 
	&+\left( v'' + \frac{2(1+ \alpha)v'}{r} - \frac{(6-\alpha(\alpha+1))v}{r^2}\right)v \eta^2 r^2 \Big]\, dr \nonumber \\
	&=\int_0^\infty \left( v^2|\eta'|^2r^2 +{2\alpha v v'}\eta^2{r} + {(4+\alpha(\alpha+1))v^2 \eta^2}\right)\, dr \nonumber \\
	&\geq \int_0^\infty \left[\frac{2 \alpha (u' r -\alpha u)}{u} + 4+ \alpha(\alpha+1) \right]w_2^2\, dr \nonumber \\
	&= \int_0^\infty \left[\frac{2 \alpha u' r }{u} + 4+ \alpha(-\alpha+1) \right]w_2^2\, dr
\end{align*}

We choose $\alpha=2$ and notice that due to properties of $u$ function $v=\frac{u}{r^2}$ is smooth and bounded. Therefore Hardy decomposition makes sense and we have 
\[
\Phi_{0,2}(0,w_2,0) \geq  \int_0^\infty \left(4\frac{u' r}{u} + 2 \right)w_2^2\, dr.
\]
Lemma is proved.
\end{proof}

\begin{lemma}\label{lemma:phi20}
For any $w_0 \in C_c^\infty(0,\infty)$ the following inequality holds
\[
\Phi_{0,2} (w_0, 0, 0) \geq \int_0^\infty \left( \frac{8 u'r}{u} - \frac{44}{9} + \frac{24u}{u'r} \right)w_0^2\, dr.
\]
\end{lemma}
\begin{proof}
We want to use a Hardy decomposition for $w_0$ similar to what was done in the previous lemma. We define $v= \frac{u'}{r^\alpha} $, where $u$ is a solution of \eqref{RS::ODE0} and $\alpha \in \RR$ to be specified later.  Let us derive an equation for $v$. Differentiating \eqref{RS::ODE0} we obtain
$$
u''' + \frac{2 u''}{r}-\frac{8u'}{r^2} +\frac{12 u}{r^3} =\hat f(u)u'.
$$
It is straightforward to compute
$$
v'=\frac{u''}{r^\alpha} - \frac{\alpha u'}{r^{\alpha+1}}, \quad v'' = \frac{u'''}{r^\alpha} - \frac{2 \alpha u''}{r^{\alpha+1}}+  \frac{\alpha(\alpha+1) u'}{r^{\alpha+2}}.
$$
Combining above equalities and rearranging terms we obtain
$$
v'' + \frac{2(1+\alpha)}{r} v'  -  \frac{[8-\alpha(\alpha+1)]v}{r^2} +\frac{12 u}{r^{3+\alpha}} =  \hat f(u)v.
$$
We would like to use this equation in the Hardy decomposition $w_0=v \xi$ 
\begin{align*}
\int_0^\infty &\left( |w_0'|^2\, r^2 + \hat f(u) w_0^2 \, r^2 \right) \, dr  \\
	&=\int_0^\infty \left[ |v' \xi + v \xi'|^2 r^2 + \left( v'' + \frac{2(1+\alpha)}{r} v'  -  \frac{[8-\alpha(\alpha+1)]v}{r^2} +\frac{12 u}{r^{3+\alpha}} \right)v \xi^2 r^2 \right]\, dr \\
	&=\int_0^\infty \left( v^2|\xi'|^2r^2 + {2\alpha v v'}\xi^2{r} - {(8-\alpha(\alpha+1))v^2 \xi^2 + \frac{12 u v \xi^2}{r^{1+\alpha}} }\right)\, dr.
\end{align*}
Using the above equality we have:
\begin{align*}
\Phi_{0,2} (w_0, 0, 0) &=\int_0^\infty\left(  2 |w'_0|^2 r^2 + 24 \,|w_0|^2 + 2 \hat f(u)\,|w_0|^2 r^2 \right) \, dr \nonumber \\
	&=\int \left( 2v^2|\xi'|^2r^2 + {4 \alpha v v'}\xi^2{r} +{[8+ 2\alpha(\alpha+1)] v^2 \xi^2 + \frac{24 u v \xi^2}{r^{1+\alpha}} }\right)\, dr \nonumber \\
	&=\int \left( 2v^2|\xi'|^2r^2 + \left[ \frac{4\alpha (u'' r - \alpha u')}{u'} + [8+ 2\alpha(\alpha+1)] + \frac{24u}{u'r} \right]w_0^2 \right)\, dr \nonumber \\
	&\geq \int   \left[ \frac{4\alpha u'' r}{u'} + [8+ 2\alpha(-\alpha+1)] + \frac{24u}{u'r} \right]w_0^2 \, dr .
\end{align*}
Choosing $\alpha=\frac{2}{3}$, defining $w_2=\frac{u'}{r^{2/3}} \xi$ and using the inequality
$$
\frac{u'' r}{u'}> \frac{3 u'r}{u} -5
$$
(see relation \eqref{rel:uprimes} in Theorem~\ref{thm:ODE}), we have
\be
\Phi_{0,2} (w_0, 0, 0)
	\geq \int_0^\infty \left( \frac{8 u'r}{u} - \frac{44}{9} + \frac{24u}{u'r} \right)w_0^2\, dr\ge 0.
\ee 
Lemma is proved.
\end{proof}
\begin{lemma}\label{lemma:phi24}
For any $w_4 \in C_c^\infty(0,\infty)$ the following inequality holds
\begin{align*}
\Phi_{0,2} (0, 0, w_4) 
	&\geq \int_0^\infty\left(   - 43  + {120 u  \over u' r} + 6\left[2- \frac{u' r}{u}\right]\left[\frac{2 u' r}{u} +3\right]  \right) |w_4|^2 \, dr  \\
	&+4 \int_0^\infty\left( \tilde f(u) -\frac{5}{2} \hat f(u) +6 f(u)  \right)\,|w_4|^2 r^2 \, dr  	.
\end{align*}
\end{lemma}
\begin{proof}
We can rewrite $\Phi_{0,2}(0,0,w_4)$ as
\begin{align*}
\Phi_{0,2} (0,0,w_4)&= \int_0^\infty\left(  4 |\partial_r w_4|^2 r^2 + 16 \,|w_4|^2 + 4 \hat f(u) \, |w_4|^2 r^2  \right) \, dr \\
	&+4 \int_0^\infty\left( \frac{3}{2} \hat f(u) -6 f(u)  \right)\,|w_4|^2 r^2 \, dr \\
	&+4 \int_0^\infty\left( \tilde f(u) -\frac{5}{2} \hat f(u) +6 f(u)  \right)\,|w_4|^2 r^2 \, dr  .
\end{align*}
We use Hardy decomposition  $w_4 = u' \gamma$, where $\gamma \in C_c^\infty (0,\infty)$, and similarly as in subsection~\ref{instab} we obtain
\begin{align*}
\int_0^\infty\left(  4 |\partial_r w_4|^2 r^2 + 16 \,|w_4|^2 \right. &+ \left.4 \hat f(u) \, |w_4|^2 r^2  \right) \, dr \\
 &= \int_0^\infty\left(  4 |u'|^2|\gamma'|^2 r^2 -16 \,|u'|^2 \gamma^2 + {48 u u' \over r}  \gamma^2  \right) \, dr 
\end{align*}
From \cite{ODE_INSZ}  (see \eqref{rel:relfs} in Theorem~\ref{thm:ODE})  we know that 
$$
(\hat f(u) -3 f(u)) u' r> -2 f(u) u
$$
and  therefore
\begin{align}\label{phi2reduced}
\Phi_{0,2}(0,0,w_4)&\geq \int_0^\infty\left(  4 |u'|^2|\gamma'|^2 r^2 -16 \,|u'|^2 \gamma^2 + {48 u u' \over r}  \gamma^2  \right) \, dr \\
	&-12 \int_0^\infty f(u) u\, u'\, \gamma^2\,r\, dr - 6  \int_0^\infty f(u) \,|w_4|^2\,r^2\, dr\nonumber\\
	&+4 \int_0^\infty\left( \tilde f(u) -\frac{5}{2} \hat f(u) +6 f(u)  \right)\,|w_4|^2 r^2 \, dr. \nonumber  
\end{align}
Using equation \eqref{RS::ODE0} for $u$ we notice that 
$$
-\int_0^\infty f(u) u\, u'\, \gamma^2\,r\, dr = \int_0^\infty \frac{6}{r} u\, u'\, \gamma^2 \, dr - \int_0^\infty \frac{3}{2} |u'|^2\, \gamma^2 \, dr +  \int_0^\infty |u'|^2\, \gamma\, \gamma' \, r \, dr.
$$
We would like to combine $4 |u'|^2 |\gamma'|^2 r^2$ with $12 |u'|^2\, \gamma\, \gamma' \, r$ and notice that
$$
4|u'|^2 |\gamma'|^2 r^2 + 12 |u'|^2 \gamma \gamma' r +9 |u'|^2 |\gamma|^2 >0.
$$
Therefore, using the last two relations in \eqref{phi2reduced}, we get:
\begin{align}\label{phi2reduced+}
\Phi_{0,2} (0,0,w_4)&\geq \int_0^\infty\left(   - 43 \,|u'|^2 \gamma^2 + {120 u u' \over r}  \gamma^2  \right) \, dr - 6  \int_0^\infty f(u) \,|w_4|^2\,r^2\, dr \nonumber \\
	&+4 \int_0^\infty\left( \tilde f(u) -\frac{5}{2} \hat f(u) +6 f(u)  \right)\,|w_4|^2 r^2 \, dr  .
\end{align}
From \cite{ODE_INSZ} (see also \eqref{rel:fws} in Theorem~\ref{thm:ODE}) we know that
$$
-f(u)r^2 > \left(2- \frac{u' r}{u} \right) \left(\frac{2 u' r}{u} +3\right) \text{ for } r > 0.
$$
Plugging this in \eqref{phi2reduced+}  we obtain the result.
\end{proof}

\vskip 0.1cm

\noindent Now we are ready to prove Proposition~\ref{prop:phi2coercivity}.

\vskip 0.1cm

\begin{proof}[Proof of Proposition~\ref{prop:phi2coercivity}] Since $C_c^\infty(0,\infty)$ is dense in $H^1((0,\infty);r^2\,dr)$, it suffices to consider $w_0, w_2, w_4 \in C_c^\infty (0,\infty)$.

Let us denote $w=\frac{u' r}{u}$. Combining lemmas \ref{lemma:phi22}, \ref{lemma:phi20} and \ref{lemma:phi24} we obtain
\begin{align*}
\Phi_{0,2} (w_0, w_2, w_4) & \geq \int_0^\infty  \left( -\frac{44}{9} + \frac{24}{w} +8 w \right) w_0^2\, dr + \int_0^\infty  \left( 4w +2 \right) w_2^2\, dr\\
	& + \int_0^\infty \left(    - 43  + {120  \over w} + 6(2-w)(2w +3)  \right) w_4^2\, dr \\
	&- \int_0^\infty 24w_0w_2\, dr + \int_0^\infty 16 w_2 w_4 \, dr\\
	&+4 \int_0^\infty\left( \tilde f(u) -\frac{5}{2} \hat f(u) +6 f(u)  \right)\,|w_4|^2 r^2 \, dr .
\end{align*}
Since $f(u) <0$ it is clear that 
$$
\tilde f(u) -\frac{5}{2} \hat f(u) +6 f(u)=-\frac{9}{2}a^2+\frac{b^2}{3}u-\frac{c^2}{3}u^2 > - \frac{11}{2} a^2.
$$
On the other hand we know that $ \frac{1}{w} > {\alpha}r^2$ and $\alpha>0$ is a constant independent of $a^2$ (see \eqref{Eq1.14} and \eqref{rel:u'u1} in Theorem~\ref{thm:ODE}).  Therefore we have 
\begin{align*}
\Phi_{0,2} (w_0, w_2, w_4) & \geq \int_0^\infty  \left( -\frac{44}{9} + \frac{24}{w} +8 w \right) w_0^2\, dr + \int_0^\infty  \left( 4w +2 \right) w_2^2\, dr\\
	& + \int_0^\infty \left(    - 43  + {119.99  \over w} + 6(2-w)(2w +3)  \right) w_4^2\, dr \\
	&- \int_0^\infty 24w_0w_2\, dr + \int_0^\infty 16 w_2 w_4 \, dr +\left(  \frac{\alpha}{100} -22 a^2 \right) \int_0^\infty |w_4|^2 r^2 \, dr  .
\end{align*}
Now we define 
\begin{align*}
\alpha_1 (w)&=  -\frac{44}{9} + \frac{24}{w} +8 w, \\
\beta_1(w) &=4w +2, \\
\gamma_1(w)&=   - 43  + {119.99  \over w} + 6(2-w)(2w +3).
\end{align*}
We have to understand if there exists $\delta_0>0$ such that
$$
\alpha_1 w_0^2 + \beta_1w_2^2 +\gamma_1w_4^2 -24 w_0 w_2 +16w_2 w_4 \geq \delta_0 \left(w_0^2+w_2^2+w_4^2\right).
$$
From  relation \eqref{rel:u'u} in Theorem~\ref{thm:ODE}  we know that  $0\leq w \leq 2$ and therefore, since $\alpha_1, \beta_1, \gamma_1 >0$ on $0 \leq w \leq 2$, we have to show that
$$
g(w)=\left((\beta_1 (w)-\delta_0) - \frac{144}{(\alpha_1 (w) -\delta_0)} \right) (\gamma_1(w)-\delta_0) -64>0
$$
for $0\leq w \leq 2$. Therefore, it is enough to show that 
$$
P(w) = g(w) w^2 (\alpha_1 (w) -\delta_0) \geq 0
$$ 
for $0 \leq w \leq 2$. Noting that $P(w)$ is a polynomial of degree six (with explicit constant coefficients) one can show that it  is positive on the interval $0\leq w\leq 2$ provided $0< \delta_0 \leq \frac{1}{1000}$.
This implies
\begin{equation}\label{eq:lowerbdphi02}
\Phi_{0,2} (w_0, w_2, w_4) \geq  \delta_0 \int_0^\infty  (w_0^2+w_2^2+w_4^2) \, dr + \left(  \frac{\alpha}{100} -22 a^2 \right) \int_0^\infty |w_4|^2 r^2 \, dr .
\end{equation}
Taking $a^2$ small enough so that $22a^2 < \frac{\alpha}{100}$ we conclude the proof. 
\end{proof}

\subsection{Proof of the main theorem}
\label{subsec:stabproof}
Now we can combine all the above results to prove our main result stated in Theorem \ref{thm:main}. From Proposition \ref{prop:phi0,2reduction}, recalling how $\Phi_{0,1}$ enters into the expression for the second variation $\mcQ(V)$  we see that there is exactly three dimensional vector space $W$, where $\mcQ(V)=0$ for all $V \in W$. Since the energy of the hedgehog is invariant under translations it follows that functions $\partial_{p_i} H(x-p)$, $i=1,2,3$ are in the kernel of $\mcQ$. Therefore the kernel of the second variation coincides with $W$ and is a span of $\{ \partial_{p_i} H(x-p) \}_{i=1}^3$. Using Lemmas \ref{DVk}, \ref{Phi3,2}, Propositions \ref{Phi012}, \ref{prop:Phikdecomp},  \ref{prop:phi0,2reduction}, \ref{prop:phi2coercivity}  we see that outside of the kernel $\mcQ(V)>0$ for $a^2$ sufficiently small. This proves the first part of the theorem. The second part of the theorem follows from Theorem \ref{thm:scalar}.
\qed

\subsection{Proof of Corollary \ref{Cor:MHLocMin}}

Fix $\Omega$ as in the statement. It suffices to show that there exists $\delta_0 > 0$ which might depend on $\Omega$ such that
\[
\mcQ(V) \geq \delta_0 \|V\|_{H^1}^2 \text{ for all } V \in H^1_0(\Omega;\mcS_0).
\]

Arguing indirectly, assume that the above inequality does not hold. Then we can find a sequence $V_m \in H^1_0(\Omega;\mcS_0)$ with $\|V_m\|_{H^1} = 1$ such that $\mcQ(V_m) \rightarrow 0$. We can further assume that $V_m$ converges weakly in $H^1$ and strongly in $L^4$ to some $V_\infty$. Then
\[
\mcQ(V_\infty) \leq \liminf_{m \rightarrow \infty} \mcQ(V_\infty) = 0.
\]
By Theorem \ref{thm:main}, this implies that the extension by zero $\bar V_\infty$ of $V_\infty$ belongs to the span of $\{\partial_{x_i} H\}$, i.e. $\bar V_\infty = \alpha \cdot \nabla H$ for some vector $\alpha \in \RR^3$.

A direct computation gives
\[
\alpha \cdot \nabla H = \Big(u' - \frac{2u}{r}\Big)\frac{\alpha \cdot x}{|x|^3} x \otimes x - u' \frac{\alpha \cdot x}{3|x|}\,\Id + \frac{u}{|x|^2} [\alpha \otimes x + x \otimes \alpha],
\]
and
$$
|\alpha \cdot \nabla H|^2 =  \frac{2(\alpha \cdot x)^2}{3|x|^2} |u'|^2 +  2 u^2 \frac{|\alpha|^2|x|^2 - (\alpha \cdot x)^2 }{|x|^4} \geq 0.
$$
It is clear that $|\alpha \cdot \nabla H|^2 =0$ if and only if $\alpha =0$.
Since $\alpha \cdot \nabla H = \bar V_\infty$ vanishes outside of $\Omega$, this implies that $\alpha = 0$, whence $V_\infty \equiv 0$. Recalling that $\|V_m\|_{H^1} = 1$ and $V_m \rightarrow V_\infty \equiv 0$ strongly in $L^4$, we then deduce that $\mcQ(V_m) \rightarrow 1$, which contradicts the definition of $V_m$. The proof is complete.
\qed

\subsection{Proof of Proposition \ref{Cor:RestrictiveStability}}

First, observe that if $Q_t$ is a smooth curve of uniaxial tensors with $Q_0 = H$, then $\frac{d}{dt}Q$ is generated by $E_0$, $E_1$, $E_2$ (see section 2.1). Therefore, in view of Lemmas \ref{DVk}, \ref{Phi3,2}, Propositions \ref{Phi012}, \ref{prop:Phikdecomp},  \ref{prop:phi0,2reduction}, it suffices to show that
\begin{align*}
0 \leq \Phi_{0,2}(w_0, w_2, 0)
	& = \int_0^\infty    \Big\{2|\partial_r w_0|^2   + |\partial_r w_2|^2
			 + \frac{1}{r^2}\Big[24 |w_0|^2 + 10 |w_2|^2
					   - 24\,w_0\,w_2 \Big] \\
  			&\qquad + 2\hat f(u)\, |w_0|^2  + f(u)  |w_2|^2  \Big\}\,r^2 \,dr.
\end{align*}
where $w_0$ and $w_2$ are functions of $r$ and belongs to $C_c^\infty(0,\infty)$. To this end, we proceed as in subsection \ref{subsec:phi0,i}. Using Hardy decomposition by taking $w_0=u' \xi$ and $w_2=u \eta$, where $\xi, \eta \in C_c^\infty (0,\infty)$, we obtain
\begin{align*}
\Phi_{0,2}(w_0, w_2, 0)= \mcE_{0,2} (\xi,\eta) 
	& = \int_0^\infty    \Big\{2|u'|^2|\xi'|^2   + u^2\,|\eta'|^2\\
			&\qquad + \frac{1}{r^2}\Big[8|u'|^2\,\xi^2 + \frac{24}{r}\,u\,u'\,\xi^2+ 4u^2\,\eta^2
					   - 24\,u\,u'\,\xi\,\eta \Big]  \Big\}\,r^2 \,dr
\end{align*}
Using inequalities
$$
8|u'|^2 \xi^2 + 2 u^2 \eta^2 -8 u u' \xi \eta \geq 0,
$$
and 
$$
\frac{24}{r} u u' \xi^2 + \frac{8}{3} r u u' \eta^2 - 16 u u' \xi \eta \geq 0,
$$
we obtain
\begin{align*}
 \mcE_{0,2} (\xi,\eta) 
	& \geq \int_0^\infty    \Big\{ u^2\,|\eta'|^2
			 + \frac{1}{r^2}\Big[2u^2\,\eta^2 
					   - \frac{8}{3}r\,u\,u'\,\eta^2 \Big]  \Big\}\,r^2 \,dr\\
	& = \int_0^\infty    \Big\{ u^2\,|\eta'|^2
			 + \frac{1}{r^2}\Big[\frac{10}{3}u^2\,\eta^2 
					   + \frac{8}{3}r u^2 \eta\eta' \Big]  \Big\}\,r^2 \,dr \geq 0.			
\end{align*}
Therefore $\Phi_{0,2}(w_0, w_2, 0) \geq 0$ with equality holding if and only if $w_0=w_2=0$. Corollary is proved.
\qed

\section{Proof of Proposition \ref{prop:Phikdecomp}}
\label{sec:bases}

\medskip
\noindent\underline{\it Case 1: $k = 0$.}
Observe that \eqref{EVP:k0U2} is a particular case of Laplace's spherical harmonics (see for instance \cite[Section V.8]{Courant-Hilbert}). Therefore,
the eigenvalues of \eqref{EVP:k0U2} are simple and given by $\lambda_{0,i}= i(i+1)$ for $i\geq 1$ and the eigenfunctions $\{ u^{(2)}_{0,i} \}_{i=1}^\infty$  (which are the associated Legendre polynomials of first order $P^1_i(\cos \theta)$) form an orthonormal basis in $L^2((0,\pi);\sin\theta\,d\theta)$; see for instance \cite[Section V.10.2]{Courant-Hilbert}). By \eqref{eigenvalu}, $u^{(0)}_{0,i}$ and $u^{(4)}_{0,i}$ satisfy the following spectral problems associated to the operators $T_0^{(0)}$ and $T_0^{(4)}$ in \eqref{EVP:U0} and \eqref{EVP:U4}:
\begin{equation}
\left\{\begin{array}{l}
\displaystyle -\frac{1}{\sin\theta}\partial_\theta(\sin\theta \partial_\theta u^{(0)}_{0,i})  = \lambda_{0,i}\,u^{(0)}_{0,i}
	,\\
\displaystyle\partial_\theta u^{(0)}_{0,i}(0) = \partial_\theta u^{(0)}_{0,i}(\pi) = 0, \qquad \int_0^\pi |u^{(0)}_{0,i}|^2\,\sin\theta\,d\theta = \lambda_{0,i}
	,
\end{array}\right.
	\label{EVP:k0U0}
\end{equation}
\begin{equation}
\left\{\begin{array}{l}
\displaystyle -\frac{1}{\sin\theta}\partial_\theta(\sin\theta \partial_\theta u^{(4)}_{0,i}) +  4\,\csc^2\theta\,u^{(4)}_{0,i} = \lambda_{0,i}\,u^{(4)}_{0,i}
	,\\
\displaystyle u^{(4)}_{0,i}(0) = u^{(4)}_{0,i}(\pi) = 0, \qquad \int_0^\pi |u^{(4)}_{0,i}|^2\,\sin\theta\,d\theta = \lambda_{0,i} - 2.
\end{array}\right.
	\label{EVP:k0U4}
\end{equation}
The normalizations of $u^{(0)}_{0,i}$ and $u^{(4)}_{0,i}$ are implied by the normalization of $u^{(2)}_{0,i}$ and the identity (which is derived from \eqref{EVP:k0U2}):
\begin{equation}
\int_0^\pi\Big[ |\partial_\theta u^{(2)}_{0,i}|^2\sin \theta+\frac{|u^{(2)}_{0,i}|^2}{\sin \theta}\Big]\, d\theta=\lambda_{0, i}, \quad i\geq 1.
	\label{Eq:17VI13.E1}
\end{equation}
The Dirichlet boundary conditions satisfied by $u^{(4)}_{0,i}$ follow from the boundary conditions of $u^{(2)}_{0,i}$. 
The Neumann boundary conditions for $u^{(0)}_{0,i}$ follow the identity $\partial_\theta u^{(0)}_{0,i}=-\lambda_{0,i} u^{(2)}_{0,i}$, which is a consequence of \eqref{Def:U0} and \eqref{EVP:k0U2}.

Since $T^{(0)}_0$ and $T^{(4)}_0$ are symmetric, the eigenfunction sets $\{u^{(0)}_{0,i}\}_{i = 1}^\infty$ and $\{u^{(4)}_{0,i}\}_{i = 2}^\infty$ are orthogonal in $L^2((0,\pi);\sin\theta\,d\theta)$. (Note that since $u^{(2)}_{0,1}(\theta) =P^1_1(\cos \theta)= \frac{\sqrt{3}}{2}\sin\theta$, $u^{(4)}_{0,1} \equiv 0$ and therefore is excluded from the above set of eigenfunctions of $T^{(4)}_0$). Also, thanks to the normalizations in \eqref{EVP:k0U0} and \eqref{EVP:k0U4}, we have $u^{(0)}_{0,i} \not\equiv 0$ for $i \geq 1$ and $u^{(4)}_{0,i} \not\equiv 0$ for $i \geq 2$.

Observe that $u^{(0)}_{0,0} := 1$ is an eigenfunction associated to the eigenvalue $\lambda_{0,0}=0$ for the problem \eqref{EVP:k0U0}.  Integrating \eqref{EVP:k0U0}, one gets that $u^{(0)}_{0,0}$ is orthogonal to $ u^{(0)}_{0,i}$ for all $i \geq 1$.  
We now want to show that  $\{u^{(0)}_{0,i}\}_{i = 1}^\infty \cup \{u^{(0)}_{0,0} \equiv 1\}$
forms an orthogonal basis in $L^2((0,\pi);\sin\theta\,d\theta)$. Since the inverse of $T_0^{(0)}$ is compact (as an operator on $L^2((0,\pi);\sin\theta\,d\theta)$), it suffices to show that, up to normalization, these exhaust all eigenfunctions associated with the problem \eqref{EVP:k0U0}. Equivalently, we need to show that any solution $\psi$ of \eqref{EVP:k0U0} which is associated to an eigenvalue $\lambda$ and is orthogonal to $u^{(0)}_{0,0}$, i.e. $\int_0^\pi \psi\,\sin\theta\,d\theta =0$, gives rise through \eqref{Def:U0} to a solution of \eqref{EVP:k0U2} with the same eigenvalue.

Indeed, integrating \eqref{Def:U0}, we obtain the function
\[
\phi(\theta) = \frac{1}{\sin\theta}\int_0^\theta \psi(\tau)\,\sin\tau\,d\tau = -\frac{1}{\sin\theta}\int_\theta^\pi \psi(\tau)\,\sin\tau\,d\tau,
\]
where in the last identity we used $\int_0^\pi \psi\,\sin\theta\,d\theta = 0$. It is readily seen that $\phi(0) = \phi(\pi) = 0$. On the other hand, from \eqref{EVP:k0U0} (taking $\lambda=\lambda_{0,i}$) we have:
\[
\partial_\theta \psi = -\frac{\lambda}{\sin\theta} \int_0^\theta \psi(\tau)\,\sin\tau\,d\tau = - \lambda\,\phi.
\]
A direct computation shows that $-\frac{1}{\sin\theta}\partial_\theta(\sin\theta \partial_\theta \phi) + \csc^2\theta \phi = \lambda \phi$. We have thus shown that $\phi$ is a solution of \eqref{EVP:k0U2}. As mentioned above, this proves that $\{u^{(0)}_{0,i}\}_{i = 1}^\infty \cup \{u^{(0)}_{0,0} \equiv 1\}$ is a basis for $L^2((0,\pi);\sin\theta\,d\theta)$.

We next show that $\{u^{(4)}_{0,i}\}_{i = 2}^\infty$ forms a basis for $L^2((0,\pi);\sin\theta\,d\theta)$. We need to show that for every eigenfunction $\chi$ of \eqref{EVP:k0U4} associated to an eigenvalue $\lambda > 2$, there is a function $\phi$ which is related to $\chi$ by \eqref{Def:U4} and solves \eqref{EVP:k0U2} (with the same eigenvalue). The argument goes as before, with the exception that the formula that determines $\phi$ from $\chi$ is somewhat trickier to obtain.

Inverting \eqref{Def:U4} using $\phi(0) = 0$ we obtain
\[
\phi(\theta) = - \sin\theta \int_0^\theta \frac{\chi(\tau)}{\sin\tau}\,d\tau + C_1\sin\theta
\]
where $C_1$ is some unknown constant which we need to determine. It is clear that $\phi(\pi) = 0$.

Since $\chi$ solves \eqref{EVP:k0U4}, $\theta = 0$ and $\theta = \pi$ are zeroes of second order of $\chi$. We thus deduce that $C_1 = \partial_\theta \phi(0)$. To determine $C_1$ we use power series expansion near $\theta = 0$. Writing $\chi(\theta)=C_2\theta^2+\hot$ and expanding the previous equation in power series near $0$ we obtain: $\phi(\theta)=C_1\theta-\left(C_1/6+C_2/2\right)\theta^3+\hot$ Since we would like $\phi$ to be a solution of \eqref{EVP:k0U2}, we need to require that
\[
8\left(C_1/6+C_2/2\right) = C_1(\lambda - 2/3), \text{ i.e. } C_1 = \frac{2\partial_\theta^2 \chi(0)}{\lambda - 2}.
\]
We thus set
\[
\phi(\theta) = - \sin\theta \int_0^\theta \frac{\chi(\tau)}{\sin\tau}\,d\tau + \frac{2\partial_\theta^2 \chi(0)}{\lambda - 2}\,\sin\theta.
\]
To check that $\phi$ satisfies \eqref{EVP:k0U2}, as in the previous case, it is more desirable to express $\phi$ using a differential relation rather than an integral relation. To this end, we note that, by \eqref{EVP:k0U4},
\[
-\sin\theta \partial_\theta(\csc\theta \partial_\theta \chi) - 2\cot\theta\,\partial_\theta \chi + 4\,\csc^2\theta\,\chi = \lambda\,\chi,
\]
from which it follows that
\begin{align*}
\partial_\theta \chi
	&= - 2\sin\theta \int_0^\theta\Big[\partial_\tau( \cot\tau\,\csc\tau \, \chi(\tau)) -  \csc\tau\,\chi(\tau)\Big]\,d\tau\\
			&\qquad  - \lambda \,\sin\theta \int_0^\theta \frac{\chi(\tau)}{\sin\tau}\,d\tau + \partial_\theta^2 \chi(0)\,\sin\theta\\
	&= -2\cot\theta\,\chi - (\lambda - 2)\sin\theta \int_0^\theta \frac{\chi(\tau)}{\sin\tau}\,d\tau   + 2\partial_\theta^2 \chi(0)\,\sin\theta
	.
\end{align*}
This leads to
\[
\phi = \frac{1}{\lambda - 2}(\partial_\theta \chi + 2\cot\theta \chi).
\]
We can now check that $\phi$ satisfies \eqref{EVP:k0U2} as desired.

To make the indices uniform, we introduce the notation that $u^{(2)}_{0,0} = u^{(4)}_{0,0} = 0$. Using the above bases, we write
\[
\u_m(r,\theta) = \sum_{i=0}^\infty w^{(m)}_{0,i}(r)\,u^{(m)}_{0,i}(\theta).
\]
Plugging this expression for $\u_m(r,\theta)$ into $\Phi_0$ given in \eqref{Phidef}, using \eqref{Eq:17VI13.E1} and the analogous relations (which are derived from \eqref{EVP:k0U0} and \eqref{EVP:k0U4})
\begin{align*}
\int_0^\pi |\partial_\theta u_{0,i}^{(0)}|^2\,\sin\theta\,d\theta 
	&= \lambda_{0,i}^2,\\
\int_0^\pi \Big[|\partial_\theta u^{(4)}_{0,i}|^2\,\sin\theta +  \frac{4\,|u^{(4)}_{0,i}|^2}{\sin\theta}\Big]\,d\theta &= \lambda_{0,i}(\lambda_{0,i} - 2),
\end{align*}
and the identities \eqref{Def:U0} and \eqref{Def:U4} (with $k = 0$), we arrive at the desired decomposition
$$\Phi_0 (\u_0, \u_2, \u_4)=\sum_{i\geq 0} \Phi_0 (w^{(0)}_{0,i}\,u^{(0)}_{0,i}, w^{(2)}_{0,i}\,u^{(2)}_{0,i}, 
w^{(4)}_{0,i}\,u^{(4)}_{0,i}) = \pi \sum_{i = 0}^\infty \Phi_{0,i} (w^{(0)}_{0,i}, w^{(2)}_{0,i}, 
w^{(4)}_{0,i}).$$

\medskip
\noindent\underline{\it Case 2: $k = 1$.} The argument is analogous to the previous case. We will only indicate the key steps.

First, $u^{(0)}_{0,i}$ and $u^{(4)}_{0,i}$ satisfy the following spectral problems associated to $T_1^{(0)}$ and $T_1^{(4)}$:
\begin{equation}
\left\{\begin{array}{l}
\displaystyle -\frac{1}{\sin\theta}\partial_\theta(\sin\theta \partial_\theta u^{(0)}_{1,i}) + \csc^2\theta\,u^{(0)}_{1,i} = \lambda_{1,i}\,u^{(0)}_{1,i}
	,\\
\displaystyle u^{(0)}_{1,i}(0) = u^{(0)}_{1,i}(\pi) = 0, \qquad \int_0^\pi |u^{(0)}_{1,i}|^2\,\sin\theta\,d\theta = \lambda_{1,i}
	,
\end{array}\right.
	\label{EVP:k1U0}
\end{equation}
\begin{equation}
\left\{\begin{array}{l}
\displaystyle -\frac{1}{\sin\theta}\partial_\theta(\sin\theta \partial_\theta u^{(4)}_{1,i}) +  (4 + (2\cot\theta + \csc\theta)^2)\,u^{(4)}_{1,i} = \lambda_{1,i}\,u^{(4)}_{1,i}
	,\\
\displaystyle u^{(4)}_{1,i}(0) = u^{(4)}_{1,i}(\pi) = 0, \qquad \int_0^\pi |u^{(4)}_{1,i}|^2\,\sin\theta\,d\theta = \lambda_{1,i} - 2
	.
\end{array}\right.
	\label{EVP:k1U4}
\end{equation}

From \eqref{EVP:k1U0} (and noting again that \eqref{EVP:k1U0} is the same as \eqref{EVP:k0U2} so it is as before a particular case of Laplace's spherical harmonics), we see that $\lambda_{1,i} \in \{j(j+1): j = 1, \ldots\}$. We will see later that $\{\lambda_{1,i}\}$ exhausts all eigenvalues of $T_1^{(0)}$ and so we must have $\lambda_{1,i} = i(i+1)$.

As before $\{u^{(0)}_{1,i}\}_{i = 1}^\infty$ and $\{u^{(4)}_{1,i}\}_{i = 2}^\infty$ are orthogonal sets of $L^2((0,\pi);\sin\theta\,d\theta)$. (Note that $u^{(2)}_{1,1}(\theta) = \frac{\sqrt{6}}{4}(1 - \cos\theta)$, and so $u^{(4)}_{1,1}$ is identically zero and hence excluded.) To show that they are bases, we need to show that, up to normalization, they exhaust all eigenfunctions associated to the problems \eqref{EVP:k1U0} and \eqref{EVP:k1U4}.

Let us start by showing that any solution $\psi$ of \eqref{EVP:k1U0} associated to some eigenvalue $\lambda$ gives rise, by means of \eqref{Def:U0}, to a solution $\phi$ of \eqref{EVP:k1U2}. Indeed, inverting \eqref{Def:U0} under the assumption that $\phi$ is regular at $\theta = 0$ we get
\[
\phi(\theta) = \frac{1}{1 - \cos\theta}\int_0^\theta (1-\cos\tau)\,\psi(\tau)\,d\tau.
\]
Since $\psi(0) = \psi(\pi) = 0$, it follows that $\partial_\theta\phi(0) = \partial_\theta\phi(\pi) = 0$. To show that $\phi$ solves \eqref{EVP:k1U2}, we use \eqref{EVP:k1U0} to find a representation of $\phi$ using differential relations rather than integral relations. By \eqref{EVP:k1U0},
\[
\frac{1}{1 - \cos\theta} \partial_\theta ((1 - \cos\theta)\partial_\theta \psi)
	= \csc\theta\,\partial_\theta \psi +\csc^2\theta\, \psi - \lambda\,\psi,
\]
which implies
\begin{align*}
(1-\cos\theta)\partial_\theta \psi
	&= \int_0^\theta \partial_\tau [(\csc\tau - \cot\tau) \psi(\tau)]\,d\tau 
			  - \lambda\int_0^\theta (1 - \cos\tau)\,\psi(\tau)\,d\tau\\
	&=  (\csc\theta - \cot\theta)  \psi 
			 - \lambda\int_0^\theta (1 - \cos\tau)\,\psi(\tau)\,d\tau
	{,} 
\end{align*}
and so
\[
\phi = \frac{1}{\lambda}(-\partial_\theta \psi + \csc\theta \, \psi) . 
\]
It follows that $\phi$ solves \eqref{EVP:k1U2}.

Similarly, if $\chi$ is a solution of \eqref{EVP:k1U4} associated to an eigenvalue $\lambda > 2$, we can find the compatible $\phi$ as follows. First, inverting \eqref{Def:U4} leads to
\[
\phi(\theta)
	= -(1 - \cos\theta)\int_0^\theta \frac{\chi(\tau)}{1 - \cos\tau}\,d\tau + C_3(1 - \cos\theta).
\]
In order that $\phi$ solves \eqref{EVP:k1U2} formally near $\theta = 0$, we need $C_3 = \frac{2\partial_\theta^3 \chi(0)}{\lambda - 2}$. To derive a differential relation, we note that \eqref{EVP:k1U4} implies
\begin{align*}
(1 - \cos\theta)\partial_\theta \Big(\frac{\partial_\theta \chi}{1 - \cos\theta} \Big)
	= - (2\cot\theta + \csc\theta )\partial_\theta \chi +  (2\cot\theta + \csc\theta)^2\,\chi - (\lambda - 4)\chi.
\end{align*} which can be integrated to obtain:
\begin{align*}
\frac{\partial_\theta \chi}{1 - \cos\theta} 
	&= - \int_0^\theta \Big[\partial_\tau \Big(\frac{2\cot\tau + \csc\tau}{1 - \cos\tau}  \chi(\tau)\Big) +  \frac{ 2}{1 - \cos\tau}\,\chi(\tau)\Big]\,d\tau\\
				&\qquad - (\lambda - 4)\int_0^\theta \frac{\chi(\tau)}{1 - \cos\tau}\,d\tau 
						+ \partial_\theta^3 \chi(0)\\
	&= - \frac{2\cot\theta + \csc\theta}{1 - \cos\theta}\,\chi 
				 - (\lambda - 2)\int_0^\theta \frac{\chi(\tau)}{1 - \cos\tau}\,d\tau		
						+ 2 \partial_\theta^3 \chi(0)	
						. 
\end{align*}
We then arrive at
\[
\phi = \frac{1}{\lambda - 2}(\partial_\theta \chi + (2\cot\theta + \csc\theta)\chi),
\]
which can be used to check that $\phi$ satisfies \eqref{EVP:k1U2}.

The decomposition for $\Phi_1$ also follows.

\medskip
\noindent\underline{Case 3: $k = 2$.} 

As in the previous cases, we start by observing that $u^{(0)}_{2,i}$ and $u^{(4)}_{2,i}$ satisfy
\begin{equation}
\left\{\begin{array}{l}
\displaystyle -\frac{1}{\sin\theta}\partial_\theta(\sin\theta \partial_\theta u^{(0)}_{2,i}) + 4\csc^2\theta\,u^{(0)}_{2,i} = \lambda_{2,i}\,u^{(0)}_{2,i}
	,\\
\displaystyle  u^{(0)}_{2,i}(0) =  u^{(0)}_{2,i}(\pi) = 0, \qquad \int_0^\pi |u^{(0)}_{2,i}|^2\,\sin\theta\,d\theta = \lambda_{2,i}
	,
\end{array}\right.
	\label{EVP:k2U0}
\end{equation}
\begin{equation}
\left\{\begin{array}{l}
\displaystyle -\frac{1}{\sin\theta}\partial_\theta(\sin\theta \partial_\theta u^{(4)}_{2,i}) +  (4 + (2\cot\theta + 2\csc\theta)^2)\,u^{(4)}_{2,i} = \lambda_{2,i}\,u^{(4)}_{2,i}
	,\\
\displaystyle u^{(4)}_{2,i}(0) = u^{(4)}_{2,i}(\pi) = 0, \qquad \int_0^\pi |u^{(4)}_{2,i}|^2\,\sin\theta\,d\theta = \lambda_{2,i} - 2
	.
\end{array}\right.
	\label{EVP:k2U4}
\end{equation}

Comparing \eqref{EVP:k2U0} and \eqref{EVP:k0U4}, we see that $\lambda_{2,i} \in \{j(j+1): j \geq 2\}$ (since $\lambda_{0,1} = 2$ gives $u_{0,1}^{(4)} \equiv 0$). That $\lambda_{2,i} = (i+1)(i+2)$ follows from the fact that $\{\lambda_{2,i}\}$ exhausts all eigenvalues associated to \eqref{EVP:k2U4}, which we will show below.

The proof then proceeds in exactly the same way as in the previous two cases. We will only show the inversion process.

Assume that $\psi$ is a solution of \eqref{EVP:k2U0} with eigenvalue $\lambda$ and $\phi$ is related to $\psi$ through \eqref{Def:U0}. Inverting \eqref{Def:U0} gives
\[
\phi(\theta)
	= \frac{\sin\theta}{(1 - \cos\theta)^2} \int_0^\theta \frac{(1 - \cos\tau)^2}{\sin\tau}\,\psi(\tau)\,d\tau.
\]
On the other hand \eqref{EVP:k2U0} gives
\[
\frac{\sin\theta}{(1 - \cos\theta)^2}\partial_\theta \Big(\frac{(1 - \cos\theta)^2}{\sin\theta}\partial_\theta \psi\Big) = 2\csc\theta\, \partial_\theta \psi + 4\csc^2\theta \,\psi - \lambda\,\psi,
\]
which implies
\begin{align*}
\partial_\theta \psi
	&= \frac{\sin\theta}{(1-\cos\theta)^2}\int_0^{\theta} \partial_\tau \Big(\frac{2(1 - \cos\tau)^2}{\sin^2\tau}\,  \psi(\tau)\Big)\,d\tau\\
			&\qquad - \lambda\,\frac{\sin\theta}{(1-\cos\theta)^2}\int_0^\theta \frac{(1 - \cos \tau)^2}{\sin\tau}\,\psi(\tau)\,d\tau\\
	&= 2\csc\theta \,  \psi(\theta)
			  - \lambda\,\frac{\sin\theta}{(1-\cos\theta)^2}\int_0^\theta \frac{(1 - \cos \tau)^2}{\sin\tau}\,\psi(\tau)\,d\tau.\end{align*}
This leads to
\[
\phi = \frac{1}{\lambda}(-\partial_\theta \psi + 2\csc\theta\,\psi)
\]
and so $\phi$ satisfies \eqref{EVP:k2U2}.

Likewise, if $\chi$ solves \eqref{EVP:k2U4} with eigenvalue $\lambda$, then inverting \eqref{Def:U4} gives
\[
\phi(\theta)
	= - \frac{(1-\cos\theta)^2}{\sin\theta} \int_0^\theta \frac{\sin\tau}{(1 - \cos\tau)^2}\,\chi(\tau)\,d\tau + C_4 \,\frac{(1-\cos\theta)^2}{\sin\theta}
	,
\]
where $C_4 = \frac{2\partial_\theta^4 \chi(0)}{3(\lambda - 2)}$ so that \eqref{EVP:k2U2} can be formally satisfied at $\theta = 0$. To derive the differential relation that determines $\phi$ in terms of $\chi$, we note that \eqref{EVP:k2U4} implies 
\begin{align*}
\frac{(1-\cos\theta)^2}{\sin\theta}\partial_\theta \Big(\frac{\sin\theta}{(1-\cos\theta)^2}\partial_\theta \chi\Big)
	&= -2(\cot\theta+\csc\theta)\partial_\theta \chi\\
			&\qquad\qquad + 4(\cot\theta + \csc\theta)^2\,\chi - (\lambda - 4)\chi.
\end{align*}
We then have
\begin{align*}
 \partial_\theta \chi
	&= -  \frac{2(1-\cos\theta)^2}{\sin\theta}\int_0^\theta  \partial_\tau \Big[(\cot\tau + \csc\tau)^3\,\csc\tau\,  \chi(\tau) \Big]\,d\tau\\
			&\qquad - (\lambda - 2)\,\frac{(1-\cos\theta)^2}{\sin\theta} \int_0^\theta \frac{\sin\tau}{(1 - \cos\tau)^2}\,\chi(\tau)\,d\tau 
					+ \frac{2\partial_\theta^4 \chi(0)}{3} \frac{(1-\cos\theta)^2}{\sin\theta}\\
	&= -2  (\cot\theta + \csc\theta)   \chi(\theta) \\
			&\qquad - (\lambda - 2)\,\frac{(1-\cos\theta)^2}{\sin\theta}\int_0^\theta \frac{\sin\tau}{(1 - \cos\tau)^2}\,\chi(\tau)\,d\tau
					 + \frac{2\partial_\theta^4 \chi(0)}{3} \frac{(1-\cos\theta)^2}{\sin\theta}
		,
\end{align*}
and so
\[
\phi = \frac{1}{\lambda - 2}(\partial_\theta \chi + 2(\cot\theta + \csc\theta) \chi).
\]
This implies that $\phi$ satisfies \eqref{EVP:k2U2}.

We have thus completed the proof.
\eproof

\appendix
\section{A Hardy-type decomposition}
\label{sec:Hardy_decomposition}

\begin{lemma}[Hardy's decomposition]\label{thm:hardy}
Let $\omega\subset \RR^N$ be an open set, $A:\omega\to \RR^{N\times N}$ be a $C^1$ nonnegative definite symmetric form, i.e., $A(x) \xi\cdot \xi\geq 0$ for every $\xi\in \RR^N$ and $x\in \omega$ and $V\in L_{loc}^1(\omega, \RR)$. We define the operator
$$L:=-\nabla\cdot (A(x)\nabla)+V(x)$$
and consider a smooth function $\psi:\omega\to \RR$ satisfying 
$$\psi>0 \,\,  \textrm{in} \, \, \omega.$$
For every $f\in C^\infty_c(\omega)$, writing $g:=\frac f \psi\in C^\infty_c(\omega)$, we have the following Hardy's decomposition: 
$$\int_\omega Lf\cdot f \, dx=\int_\omega \psi^2 A(x) \nabla g \cdot \nabla g\, dx+\int_\omega g^2 L\psi\cdot \psi \, dx.$$
In particular, if $L\psi\geq 0   \,\,  \textrm{a.e. in } \, \, \omega$, then
$$\int_\omega Lf\cdot f \, dx\geq 0  \,\,  \textrm{for every } \, \, f\in C^\infty_c(\omega). $$
\end{lemma}

\bproof
For every $f\in C^\infty_c(\omega)$, set $g:=\frac f \psi\in C^\infty_c(\omega).$ Integrating by parts, we obtain:
\begin{align*}
\int_\omega Lf\cdot f \, dx &= \int_\omega L(\psi g)\cdot (\psi g) \, dx\\
&=\int_\omega \bigg(A(x) \nabla(\psi g)\cdot \nabla (\psi g) +V(x) \psi^2 g^2\bigg)\, dx\\
&=\int_\omega \bigg(g^2 A(x) \nabla\psi \cdot \nabla \psi+  \underbrace{\psi^2 A(x) \nabla g \cdot \nabla g}_{\geq 0} +\frac 1 2 A(x) \nabla (\psi^2) \cdot \nabla (g^2) +V(x) \psi^2 g^2\bigg)\, dx\\
&{=} \int_\omega \psi^2 A(x) \nabla g \cdot \nabla g\, dx+\int_\omega \bigg(g^2 A(x) \nabla\psi \cdot \nabla \psi-\frac 1 2 g^2\nabla\cdot \big (A(x) \nabla (\psi^2)\big) +V(x) \psi^2 g^2\bigg)\, dx\\ 
&=\int_\omega \psi^2 A(x) \nabla g \cdot \nabla g\, dx+\int_\omega g^2 L\psi\cdot \psi \, dx.
\end{align*}
\eproof

\medskip
\noindent{\bf Proof of inequality \eqref{PWH}.} When $w \in C^\infty_c(0,\pi)$, the assertion follows from Lemma \ref{thm:hardy} with $\omega = (0,\pi) \subset \RR^1$, $A = \sin\theta$, $V = \frac{k^2}{\sin\theta}$ and $\psi = \sin^k\theta$. (Note that $L\psi = (k^2 + k) \sin\theta\, \psi$.) The general case follows from standard density argument.
\eproof

\section*{Acknowledgment.} The authors gratefully acknowledge the hospitality and partial support of Hausdorff Research Institute for Mathematics, Bonn where part of this work was carried out. The hospitality of the Isaac Newton Institute, Cambridge is also acknowledged. R.I.  acknowledges partial support by the ANR project ANR-10-JCJC 0106. V.S. acknowledges support by EPSRC grant  EP/I028714/1.

\bibliographystyle{siam}
\bibliography{paris,%
LiquidCrystals}

\end{document}